\documentclass[preprint,3p,sort&compress,final,times]{elsarticle}
\usepackage{amsmath}
\usepackage{amssymb}
\usepackage{mathtools}
\usepackage{placeins}
\usepackage[usenames,dvipsnames]{color}
\usepackage{xcolor}
\usepackage{cases}
\usepackage{array}

\usepackage[usenames,dvipsnames]{color}

\begin{document}
\begin{frontmatter}

\title{Petrov-Galerkin flux upwinding for mixed mimetic spectral elements, and its application to
geophysical flow problems}

\author[MON]{David Lee\corref{cor}}
\ead{davelee2804@gmail.com}

\address[MON]{Department of Mechanical and Aerospace Engineering, Monash University, Melbourne 3800, Australia}
\cortext[cor]{Corresponding author. Tel. +61 452 262 804.}

\begin{abstract}
Upwinded mass fluxes are described and analysed for advection operators discretised using mixed 
mimetic spectral elements. This involves a Petrov-Galerkin formulation by which the mass 
flux test functions are evaluated at downstream locations along velocity characteristics.
As for the original mixed mimetic spectral element advection operator, the upwinded mass 
flux advection operator is conservative, however unlike the original advection operator, 
which is purely hyperbolic, the upwinded advection operator adds dissipation which is 
biased towards high wave numbers. The upwinded advection operator
also removes the spectral gaps present in the dispersion relation for the original advection operator. 
As for the original advection operator, a material form advection operator may be constructed by similarly
downwinding the trial functions of the tracer gradients.
Both methods allow for the recovery of exact energy conservation for an incompressible flow field 
via skew-symmetric formulations. However these skew-symmetric formulations are once again purely 
hyperbolic operators which do not suppress oscillations. The scheme is implemented within
a shallow water code on the sphere in order to diagnose and interpolate the potential vorticity. 
In the absence of other dissipation terms, it is shown to yield more coherent results for a 
standard test case of barotropic instability.
\end{abstract}


\end{frontmatter}

It is well known that in the absence of upwinding, diffusion or other stabilisation schemes, 
discrete Eulerian approximations to hyperbolic terms
result in spurious oscillations. Numerous methods have been devised to address this issue in the context of
finite element methods, including the streamwise-upwind Petrov-Galerkin method \cite{BH82} and
variational multiscale methods \cite{Hughes98}. These schemes have subsequently been applied within collocated
spectral element discretisations \cite{Marras12}. In the present work we address this issue 
in the context of a mixed mimetic spectral element discretisation \cite{Gerritsma11}, a high order finite element method for
which the $H(\mathrm{div},\Omega)$ function space of the mass flux is compatibly mapped to the $L^2(\Omega)$ 
space of the tracer field by the divergence operator \cite{Kreeft13}. 

Existing variational approaches to upwind stabilisation typically involve augmenting the test space with 
additional terms, for example a convective term that results in some form of symmetric positive definite 
operator which acts as a viscosity in the direction of the flow \cite{BH82}, or additionally a dual space in 
which the fine scales are represented via an adjoint problem from which the stabilisation terms are derived \cite{Hughes98}. 
An alternative form of upwinding for mixed finite elements for Hamiltonian systems has been introduced 
whereby all of the test and trial functions in the $H(\mathrm{div},\Omega)$ space are projected onto upwind 
variants within the skew-symmetric operator so as to preserve energy conservation \cite{WCB20}.
Here we instead base our methods around the idea of Lagrangian basis functions, which are evaluated at
upstream locations in order to stabilise the resultant mass fluxes.

Mass conserving variational schemes based on Lagrangian basis functions have been developed previously for 
both finite element \cite{RC02} and discontinuous Galerkin \cite{Guo14,Lee16,Bosler19} methods. In each of the above 
cases the computation of fluxes via the extrusion of the test functions over the tracer field implies an 
adjoint problem in which the test functions are themselves advected along Lagrangian characteristics. In 
the present mixed finite element context, only the test functions for the mass fluxes are advected, and 
mass conservation is ensured since the trial functions, which compatibly map to the trial functions for the 
tracer via the divergence operator, remain static.

In the following section the advection problem is introduced for the standard mixed mimetic spectral 
element discretisation. In Section 2 variations on this method using upwinded test functions for flux 
form advection, or downwinded trial functions for material form advection will be described. Results for 
these new formulations, and comparisons to the original scheme will be presented in Section 3. In
Section 4 the scheme is extended to two dimensions by way of some standard test cases for passive advection
on the sphere.
As a practical demonstration of the Petrov-Galerkin upwind stabilisation, results are presented for the
diagnosis and interpolation of potential vorticity for the rotating shallow water equations on the sphere in 
Section 5, and conclusions are discussed in Section 6.

\section{Advection using mixed mimetic spectral elements}

Consider the scalar advection equation of a tracer, $q$, subject to a velocity field, $\boldsymbol u$, expressed in
conservative flux form as
\begin{subequations}\label{adv_eqn}
\begin{align}
\dot{q} + \nabla\cdot\boldsymbol F &= 0,\\
\qquad\boldsymbol{F} &= \boldsymbol{u}q,
\end{align}
\end{subequations}
where $\boldsymbol{F}$ is the mass flux, 
within a periodic spatial domain $\Omega = [0,L) \subset \mathbb R$, and a temporal coordinate $t = [0,\infty)$.
We further assume the existence of two discrete, finite dimensional subspaces, 
$\mathcal U_h\subset H(\mathrm{div},\Omega)$ and $\mathcal Q_h\subset L^2(\Omega)$, such that $\mathcal Q_h$ 
contains a set of functions which are square integrable, and $\mathcal U_h$ contains a set of functions for 
which the sum of the functions and their divergence are square integrable. These spaces satisfy a 
compatibility property of the form
\begin{equation}\label{H_div_L_2}
\mathcal U_h \stackrel{\nabla\cdot}{\longrightarrow} \mathcal Q_h.
\end{equation}
These two discrete subspaces are composed of a finite set of polynomial basis functions of degree $p$ within 
each element such that
\begin{equation}
\mathcal U_h = \mathrm{span}\{l^p_0(\xi),\dots,l^p_p(\xi)\},\quad\mathcal Q_h = \mathrm{span}\{e^p_0(\xi),\dots,e^p_{p-1}(\xi)\},
\end{equation}
where $\xi$ is a local coordinate within the canonical element domain $[-1,1]\subset\mathbb R$. For the
remainder of this article, we will use as these subspaces the Lagrange polynomials of degree $p$ with their
roots as the Gauss-Lobatto-Legendre (GLL) points of equal degree, and the associated edge functions \cite{Gerritsma11}
respectively. These are given for the GLL nodes $\xi_k$ as
\begin{equation}\label{eq::nodal_edge}
l^{p}_{i}(\xi) = \prod_{\substack{k=0\\k\neq i}}^{p}\frac{\xi-\xi_{k}}{\xi_{i}-\xi_{k}},\qquad
e^{p}_{i}(\xi) = - \sum_{k=0}^{i-1}\frac{\mathrm{d}l^{p}_{k}(\xi)}{\mathrm{d}\xi}.
\end{equation}
The discrete mass flux, $\boldsymbol F_h\in\mathcal U_h$, and tracer field, $q_h\in\mathcal Q_h$ may be interpolated
via the nodal and edge bases respectively as
\begin{equation}
\boldsymbol F_h(\xi) = \sum_{i=0}^p\hat F_il^p_i(\xi),\quad q_h(\xi) = \sum_{i=0}^{p-1} \hat q_ie_i^p(\xi),
\end{equation}
where $\hat F_i$ and $\hat q_i$ are the degrees of freedom.
For the remainder of this article we will drop the superscripts, $l_i^p(\xi)$ and $e_i^p(\xi)$, and assume that
$p$ remains fixed at some specified degree for a given configuration.

Due to the orthogonality and compatibility properties of $\mathcal U_h$ and $\mathcal Q_h$, by which the fundamental
theorem of calculus is satisfied exactly between GLL (or any other choice of) nodes, the discrete 
divergence operator may be represented by a purely topological relation, known as an incidence matrix \cite{Kreeft13}, 
and defined here over a periodic one dimensional domain as
\begin{equation}
\boldsymbol{\mathsf{E}} = 
\begin{bmatrix}
-1 &  1 &  0 &  0 & \ldots & 0 \\
 0 & -1 &  1 &  0 & \ldots & 0 \\
 0 &  0 & -1 &  1 & \ldots & 0 \\
\vdots & \vdots & \vdots & \vdots & \ddots & \vdots \\
1 & 0 & 0 & 0 & \ldots & -1 \\
\end{bmatrix}
.
\end{equation}
The incidence matrix provides a strong form mapping between degrees of freedom in $\mathcal U_h$ and those in 
$\mathcal Q_h$, and as such is a discrete representation of \eqref{H_div_L_2}.

The discrete divergence of the mass flux onto a semi-discrete time derivative of the tracer field is then given in 
the strong form as
\begin{equation}\label{adv_disc_1}
\dot{\hat q}_i = -\mathsf{E}_{ij}\hat F_j.
\end{equation}
We wish to derive a mass flux that is $C^0$ continuous across element boundaries, and so compute this
via a contraction of the tracer field onto the velocity field as
\begin{equation}\label{adv_disc_2}
\langle l_i,l_j\rangle_{\Omega}\hat F_j = \langle l_i\cdot\boldsymbol u_h,e_k\rangle_{\Omega}\hat q_k,
\qquad\forall l_i\in\mathcal{U}_h,
\end{equation}
where the brackets $\langle a,b\rangle_{\Omega} = \int ab\mathrm{d}\Omega$ correspond to a bilinear operator
and the subscript $\Omega$ denotes the integration and assembly over all elements in the domain $\Omega$.
Note that the set of nodal basis functions $l_i(\xi)$ have the same representation in both local and global
coordinates, and so the transformation of these bases to global coordinates is unity, and does not involve the Jacobian,
while the transformation for $e_i(\xi)\in\mathcal Q_h$ is the inverse of the Jacobian determinant, $|J(\xi)|^{-1}$. 
In higher dimensions these are transformed via the Piola mappings \cite{RKL10,LP18}.

Equations \eqref{adv_disc_1} and \eqref{adv_disc_2} describe the semi-discrete integration of \eqref{adv_eqn} 
\cite{PG13,LPG18}. These may be expressed in a single equation as
\begin{equation}\label{adv_disc_coupled}
\langle e_i,e_j\rangle_{\Omega}\dot{\hat q}_j + 
\langle e_i,e_k\rangle_{\Omega}\mathsf{E}_{km}\langle l_n,l_m\rangle_{\Omega}^{-1}\langle l_n\cdot\boldsymbol u_h,e_r\rangle_{\Omega}\hat q_r = 0
\qquad\forall e_i\in\mathcal{Q}_h,
\end{equation}
where both sides of \eqref{adv_disc_1} have been pre-multiplied by the $\mathcal Q_h$ mass matrix, 
$\boldsymbol{\mathsf M} = \langle e_i,e_j\rangle_{\Omega}$. This equation conserves mass due to the telescopic property
of the incidence matrix, $\boldsymbol{1}^{\top}\boldsymbol{\mathsf{E}} = \boldsymbol{0}$ \cite{LPG18}.
Note that in the multi-dimensional case the discrete subspace of $H(\mathrm{div},\Omega)$ is composed of vector functions with
continuous normal components, and so $l_i$ is replaced with a vector field basis with $C^0$ continuity across
element boundaries.

Notably, \eqref{adv_disc_coupled} is the adjoint of the discrete material form of the advection equation
\begin{subequations}\label{adv_material}
\begin{align}
\dot q + \boldsymbol u\cdot\boldsymbol{G} &= 0,\\
\boldsymbol{G} &= \nabla q,
\end{align}
\end{subequations}
where $\boldsymbol{G}$ is the tracer gradient. The degrees of freedom
of the discrete form of the tracer gradient, $\hat{G}_i$, are determined via a weak form integration 
by parts relation with respect to the strong form divergence operator (assuming periodic boundary conditions) as \cite{GP16}
\begin{equation}\label{disc_grad}
\hat{G}_j = -\langle l_i,l_j\rangle_{\Omega}^{-1}\mathsf{E}_{k,i}^{\top}\langle e_k,e_l\rangle\hat{q}_l\qquad\forall l_i\in\mathcal{U}_h.
\end{equation}
The discrete variational form of \eqref{adv_material}, combined with \eqref{disc_grad} then gives the 
discrete form of the material advection equation as
\begin{equation}\label{adv_disc_material}
\langle e_i,e_j\rangle_{\Omega}\dot{\hat q}_j -
\langle e_i,\boldsymbol u_h\cdot l_k\rangle_{\Omega}\langle l_m,l_k\rangle^{-1}_{\Omega}\mathsf{E}_{nm}^{\top}\langle e_n,e_r\rangle_{\Omega}\hat q_r = 0
\qquad\forall e_i\in\mathcal{Q}_h.
\end{equation}
The adjoint property is then satisfied since the advection operator in \eqref{adv_disc_coupled}, 
\begin{equation}\label{eq::A}
\boldsymbol{\mathsf A} = \langle e_i,e_k\rangle_{\Omega}\mathsf{E}_{km}\langle l_n,l_m\rangle_{\Omega}^{-1}\langle l_n\cdot\boldsymbol u_h,e_r\rangle_{\Omega},
\end{equation}
and the corresponding operator in \eqref{adv_disc_material}, 
\begin{equation}\label{eq::B}
\boldsymbol{\mathsf B} = -\langle e_i,\boldsymbol u_h\cdot l_k\rangle_{\Omega}\langle l_m,l_k\rangle_{\Omega}^{-1}\mathsf{E}_{nm}^{\top}\langle e_n,e_r\rangle_{\Omega},
\end{equation}
are related as $\boldsymbol{\mathsf{B}} = -\boldsymbol{\mathsf{A}}^{\top}$.

As with other discretisations \cite{Morinishi10,PG17}, the adjoint property for the flux form and material form 
advection operators allows for exact energy conservation for incompressible flows 
($\boldsymbol{\mathsf{E}}\hat{\boldsymbol u}_h = \boldsymbol 0,\forall\boldsymbol{u}_h\in\mathcal{U}_h$).
This is achieved via a skew-symmetric form of the advection operator as
\begin{equation}\label{adv_ss}
\boldsymbol{\mathsf S} = \frac{1}{2}(\boldsymbol{\mathsf A} + \boldsymbol{\mathsf B}) =
\frac{1}{2}(\boldsymbol{\mathsf A} - \boldsymbol{\mathsf A}^{\top}),
\end{equation}
and a centered time integration scheme of the form
\begin{equation}\label{adv_eqn_ss}
\boldsymbol{\mathsf M}\frac{(\hat q_h^{n+1} - \hat q_h^n)}{\Delta t} + \boldsymbol{\mathsf S}\frac{(\hat q_h^{n+1} + \hat q_h^n)}{2} = 0.
\end{equation}
Pre-multiplying both sides of \eqref{adv_eqn_ss} by $(\hat q_h^{n+1} + \hat q_h^n)/2$ gives the conservation of energy as
\begin{equation}
\hat q_h^{n+1}\boldsymbol{\mathsf M}\hat q_h^{n+1} = \hat q_h^{n}\boldsymbol{\mathsf M}\hat q_h^{n},
\end{equation}
due to the skew-symmetry of $\boldsymbol{\mathsf S}$ and the bi-linearity of $\boldsymbol{\mathsf M}$.

\section{Petrov-Galerkin flux upwinding}

The advection operators described in the preceeding section have numerous appealing properties, including high order 
error convergence (for smooth solutions), mass and energy conservation, and purely imaginary eigenvalues (strict hyperbolicity).
However they are also prone to spurious oscillations in the presence of sharp, poorly
resolved gradients. In this section we describe modified formulations which smooth out these oscillations via upwinding. 
A consequence of this upwinding construction is that the eigenvalues of the operators are no longer purely hyperbolic, 
such that non-zero real eigenvalues are present that act to damp the solutions.

In order to upwind the mass flux test functions, we may equivalently evaluate these at \emph{downstream} locations
\cite{Lee16}, defined locally as
\begin{equation}\label{xi_downwind}
\xi^d = \xi + \Delta t\int_{s=0}^{1}|J(\xi(s))|^{-1}\boldsymbol u(\xi(s),s)\mathrm{d}s,
\end{equation}
where $|J(\xi(s))|^{-1}$ is the Jacobian determinant inverse and $\boldsymbol{u}(\xi(s),s)$ is the
velocity in \emph{local} element coordinates,
such that $l_i^u = l_i(\xi^d)$. 
The physical units of the time step and the Jacobian determinant are introduced in order to ensure
that the amount of upwinding is small with respect to the flow velocity. High order basis functions in $\mathcal{U}_h$
diverge in proportion to their polynomial degree outside of the canonical element domain, and so if these are
upwinded excessively then the condition number of the matrix will degrade significantly.
The choice of upwinding distance used here works well experimentally but is not unique. Just as
the traditional SUPG method \cite{BH82} is dependent on the choice of a tuning parameter, there is perhaps some
way to optimise the upwinding length scale here also. 

The mass flux is then computed using a Petrov-Galerkin formulation as
\begin{equation}
\langle l_i^u,l_j\rangle_{\Omega}\hat F_j^{PG} = \langle l_i^u\cdot\boldsymbol u_h,e_k\rangle_{\Omega}\hat q_k,
\qquad\forall l_i^u\in\mathcal{U}_h.
\end{equation}
Since the trial space remains unaltered, the resulting mass flux $\hat F_j^{PG}$ maintains its
compatible mapping with respect to the space of the tracer field, $\mathcal Q_h$, such that 
mass conservation is preserved for the upwinded flux. The choice to upwind the test functions for
the mass flux and not the full advection equation, as is customary in Petrov-Galerkin formulatons \cite{BH82,Hughes98},
is motivated by the need to upwind a test space for which $C^0$ continuity is enforced in the direction 
of the flow. If the tests function for the full advection equation \eqref{adv_disc_coupled}
are upwinded instead of those for the mass flux, then nothing is gained, since these are discontinuous
across element boundaries, and moreover the $\mathcal{Q}_h$ mass matrix $\boldsymbol{\mathsf{M}}$
can be cancelled from \eqref{adv_disc_coupled} since the divergence operator is expressed in the strong form.
The Petrov-Galerkin upwinded advection operator is then given as
\begin{equation}\label{A_PG}
\boldsymbol{\mathsf A}_{PG;\Delta t} = \langle e_i,e_k\rangle_{\Omega}\mathsf{E}_{km}\langle l_n^u,l_m\rangle_{\Omega}^{-1}\langle l_n^u\cdot\boldsymbol u_h,e_r\rangle_{\Omega}.
\end{equation}
The corresponding material form advection operator is given as
\begin{equation}\label{B_PG}
\boldsymbol{\mathsf B}_{PG;\Delta t} = -\boldsymbol{\mathsf A}_{PG;-\Delta t}^{\top} = 
-\langle e_i,\boldsymbol u_h\cdot l_k^d\rangle_{\Omega}\langle l_m,l_k^d\rangle_{\Omega}^{-1}\mathsf{E}_{nm}^{\top}\langle e_n,e_r\rangle_{\Omega}
\end{equation}
where $l_i^d = l_i(\xi^u)$ and
\begin{equation}\label{xi_upwind}
\xi^u = \xi - \Delta t\int_{s=0}^{1}|J(\xi(s))|^{-1}\boldsymbol u(\xi(s),s)\mathrm{d}s,
\end{equation}
Note that in each of the examples presented in this article we use a simple first order Euler integration in 
order to determine $\xi^d$ and $\xi^u$. For multi-dimensional flows involving large amounts of deformation, 
it may be advisable to use higher order integration in order to determine the departure locations of the Gauss-Lobatto quadrature
points. No assumption has been made about $\boldsymbol{u}_h$ in either \eqref{A_PG} or \eqref{B_PG}, other than this
be uniquely defined and that it represents a vector field for multi-dimensional domains. As such the velocity field 
may be either compressible or incompressible. The only instance in which incompressibility is required is in the 
construction of a skew-symmetric advection operator, as in \eqref{adv_ss}.

One drawback of the upwinded operators described above is that by evaluating oscillatory basis
functions outside of the canonical domain the condition number of the matrices necessary to
determine mass fluxes and tracer gradients is increased. Also for nodal bases which are orthogonal with respect 
to the standard Gauss-Lobatto quadrature points, by moving these quadrature point locations this orthogonality
is broken such that the corresponding mass matrices are no longer diagonal.

\section{Results}

\begin{figure}
\centering
\includegraphics[width=0.48\textwidth,height=0.36\textwidth]{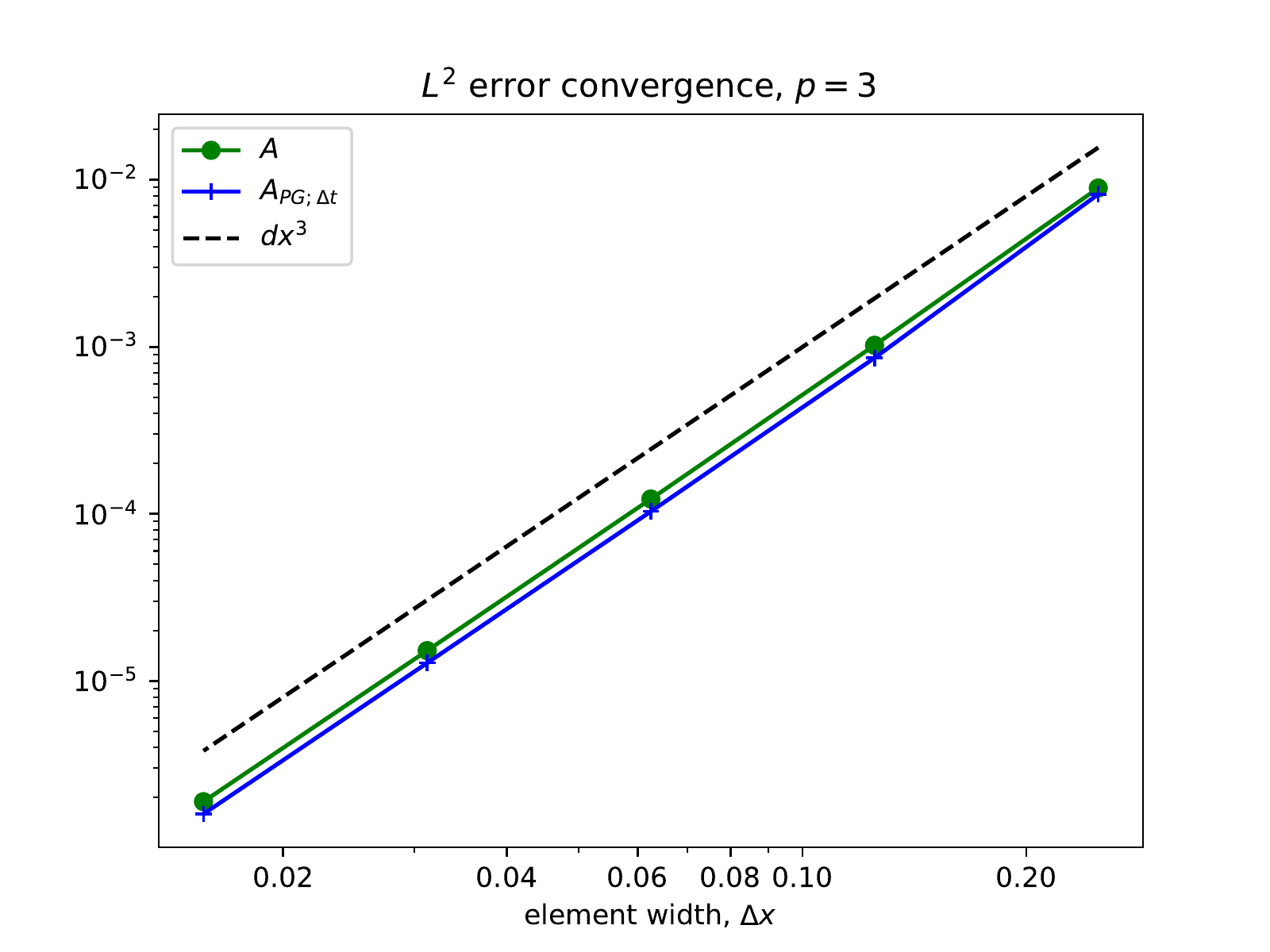}
\includegraphics[width=0.48\textwidth,height=0.36\textwidth]{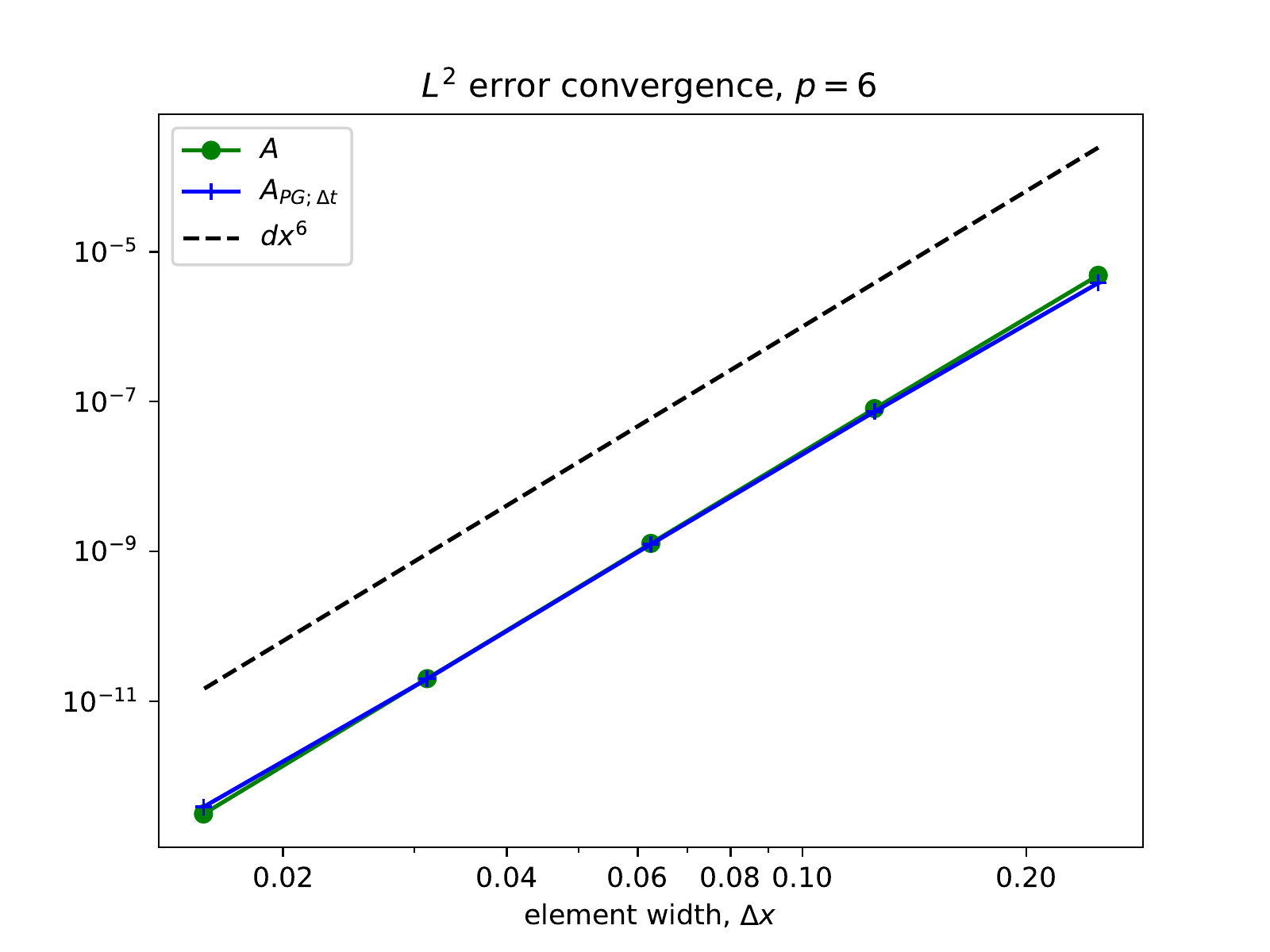}
\caption{Mass flux error convergence for $p=3$ (left) and $p=6$ (right).}
\label{fig::l2_convergence_36}
\end{figure}

We first verify the error convergence properties for the upwinded mass flux, for a manufactured solution of the form 
$q = 0.5(1.0-\cos(2\pi x))$, $\boldsymbol u = 0.4 + 0.2(1.0 + \sin(2\pi x))$ over a periodic domain of unit length $L = 1$,
such that $x\in[0,L)$, and compare the convergence rates against the original definition of the discrete mass
flux for $\boldsymbol F = \boldsymbol uq$.
For this test we compare two configurations, the first
for elements of degree $p=3$, and the second for elements of degree $p=6$. In both cases we use $n_e = 4\times 2^n$ elements,
with $n = 1,2,3,4,5$ and a time step of $\Delta t = 0.1/n_e$ in order to upwind the trial functions for the Petrov-Galerkin
formulation.
As observed in Fig. \ref{fig::l2_convergence_36}, both the original and upwinded mass fluxes converge at their theoretical
rates for both $p=3$ and $p=6$ within the $L^2(\Omega)$ norm. The convergence is marginally better for the upwinded
formulation, however this improvement diminishes with polynomial order.

\begin{figure}
\centering
\includegraphics[width=0.48\textwidth,height=0.36\textwidth]{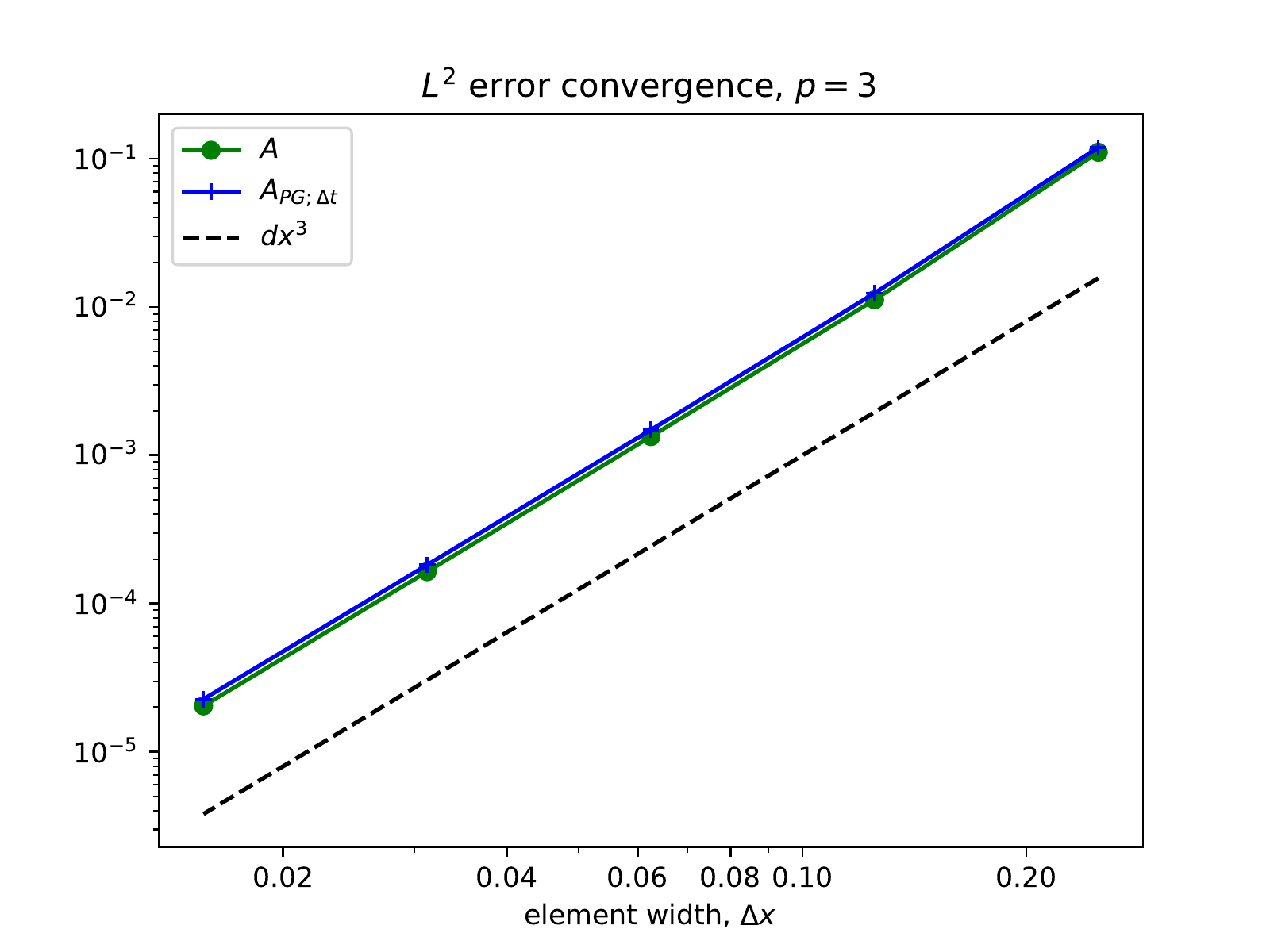}
\includegraphics[width=0.48\textwidth,height=0.36\textwidth]{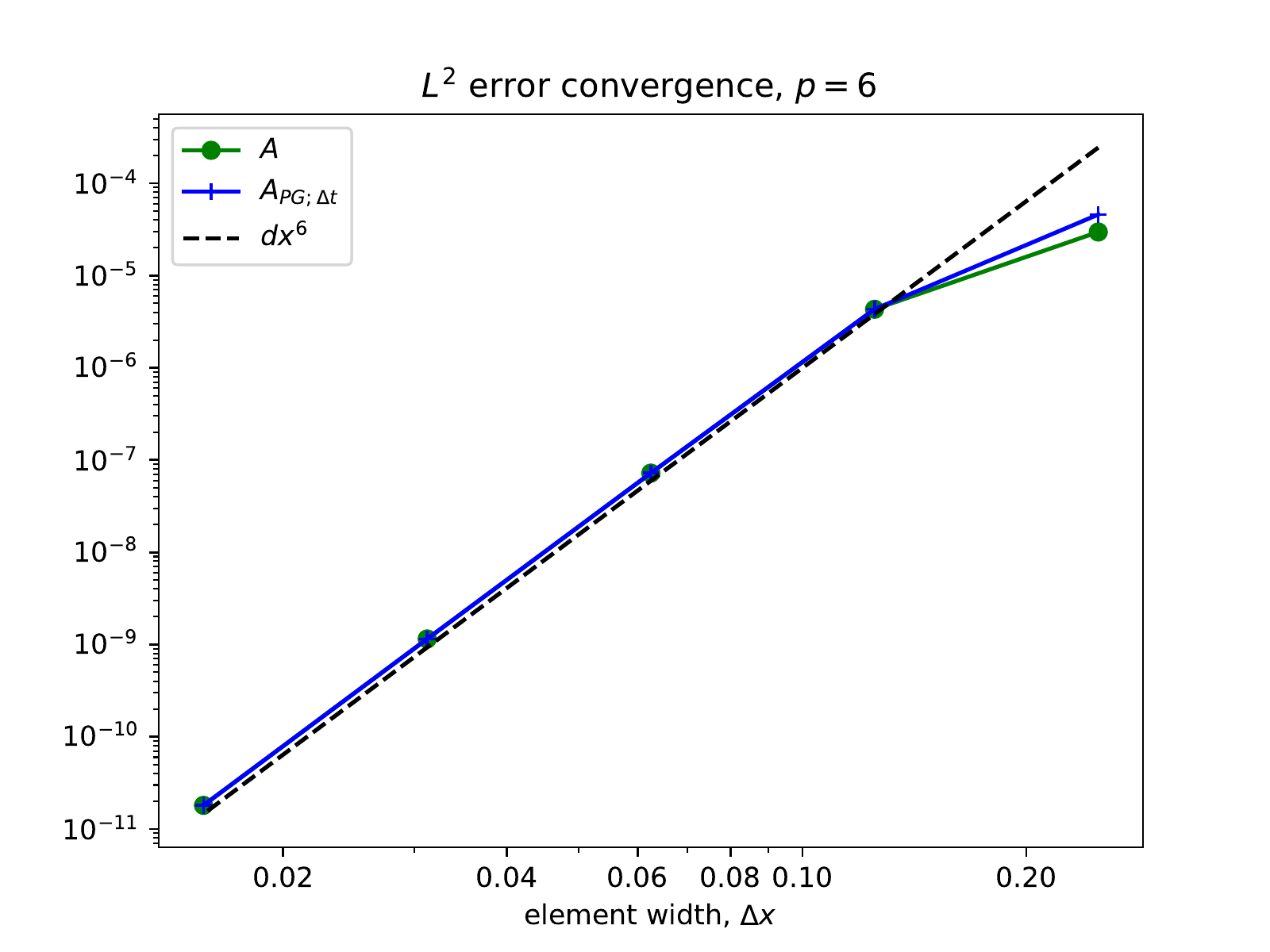}
\caption{Material advection operator error convergence for $p=3$ (left) and $p=6$ (right).}
\label{fig::l2_convergence_36_mat_form}
\end{figure}

The error convergence of the material form advection operator \eqref{eq::B} and its downwinded variant \eqref{B_PG} 
is then verified against a manufactured solution with the same specification of $q$ and $\boldsymbol{u}$ as for the
previous test case. For polynomials of degree $p=3$ and $p=6$, the errors converge at the expected rate for both 
formulations, as shown in Fig. \eqref{fig::l2_convergence_36_mat_form}.

As a second test we time step the advection equation using the original, upwinded flux form, and downwinded material form
advection operators for an incompressible flow of $\boldsymbol u = 0.4$ and an initial tracer distribution of
\begin{equation}
q(x,0) = 
\begin{cases}
0.5 + 0.5\tanh(200(x-0.4)), \qquad x < 0.5\\
0.5 + 0.5\tanh(200(0.6-x)), \qquad x\ge 0.5
\end{cases}
\end{equation}
over the unit domain of $L=1$ with 20 elements and a time step of $\Delta t = 0.005$ over a single revolution of period $T=2.5$. 
In each case a second order centered time stepping scheme of the form 
\begin{equation}\label{eq::centered_advection}
\boldsymbol{\mathsf M}\frac{(\hat q_h^{n+1} - \hat q_h^n)}{\Delta t} + \boldsymbol{\mathsf A}\frac{(\hat q_h^{n+1} + \hat q_h^n)}{2} = 0
\end{equation}
was employed (with $\boldsymbol{\mathsf A}$ replaced with $\boldsymbol{\mathsf A}_{PG;\Delta t}$ and 
$-\boldsymbol{\mathsf A}_{PG;-\Delta t}^{\top}$ for the upwinded flux form and downwinded material forms respectively).

\begin{figure}
\centering
\includegraphics[width=0.48\textwidth,height=0.36\textwidth]{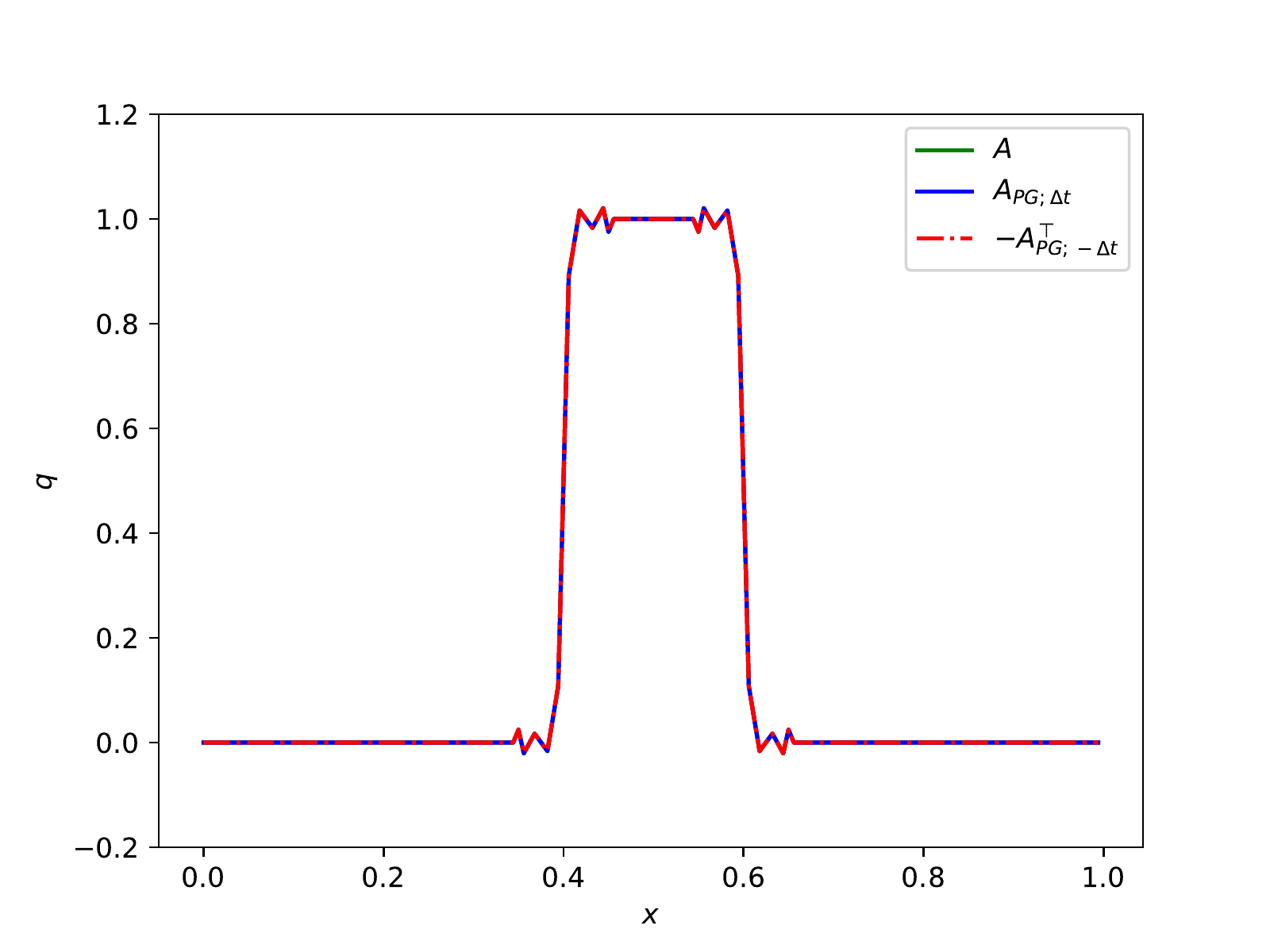}
\includegraphics[width=0.48\textwidth,height=0.36\textwidth]{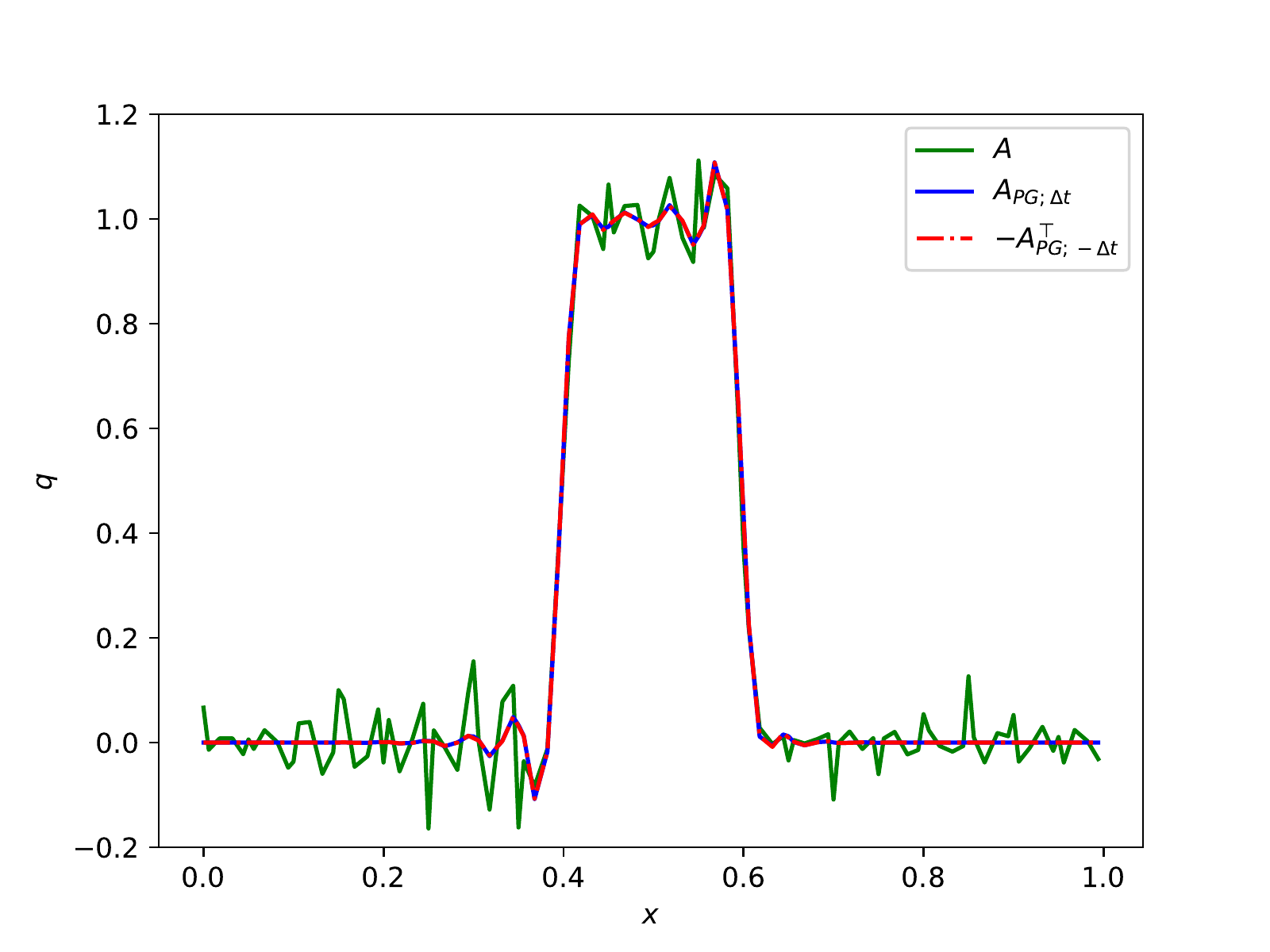}
\caption{Initial state (left) and final state (right) after one revolution, $p=5$, 20 elements, $u=0.4$, $\Delta t=0.005$.}
\label{fig::advection}
\end{figure}

\begin{figure}
\centering
\includegraphics[width=0.48\textwidth,height=0.36\textwidth]{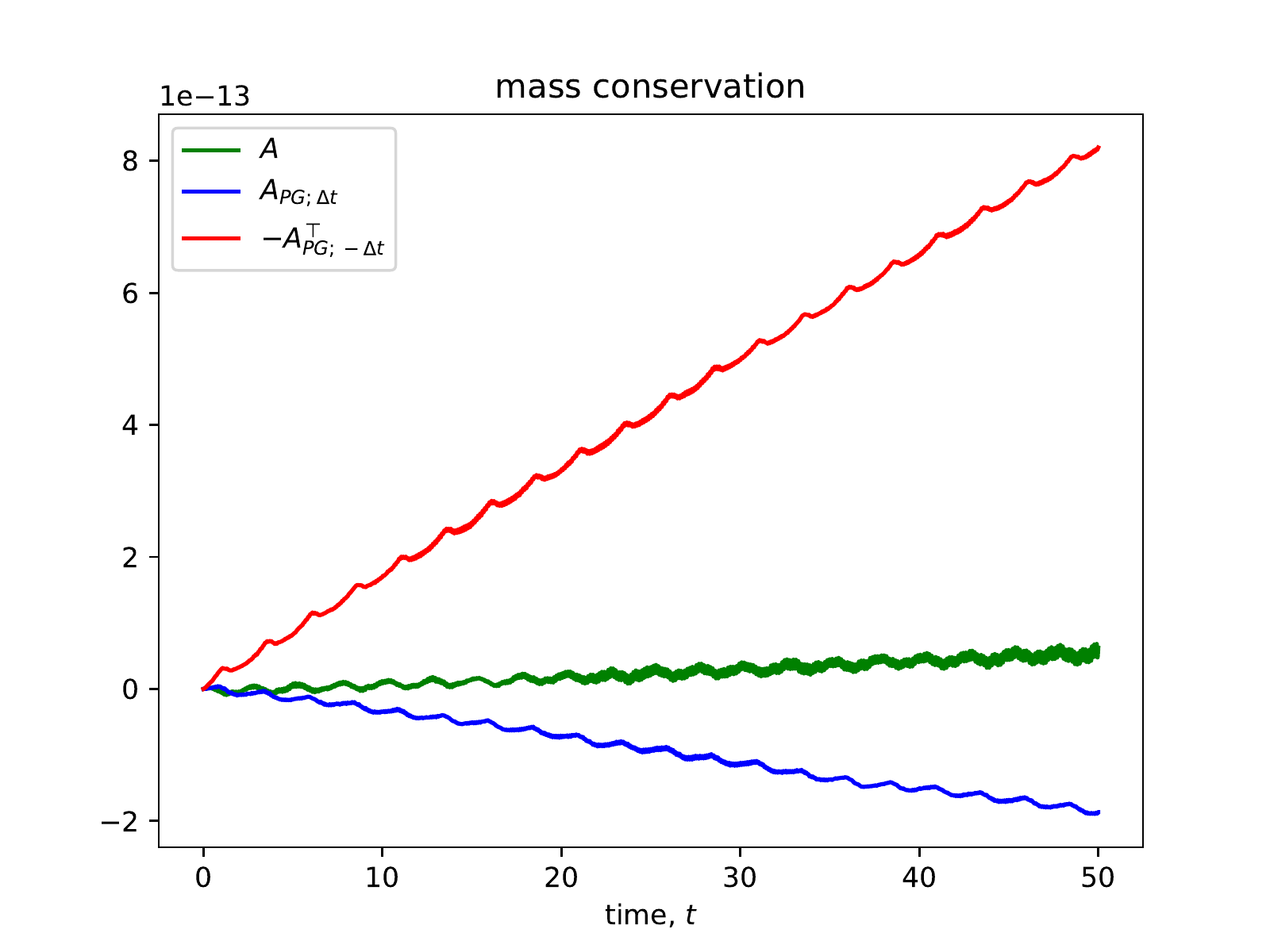}
\includegraphics[width=0.48\textwidth,height=0.36\textwidth]{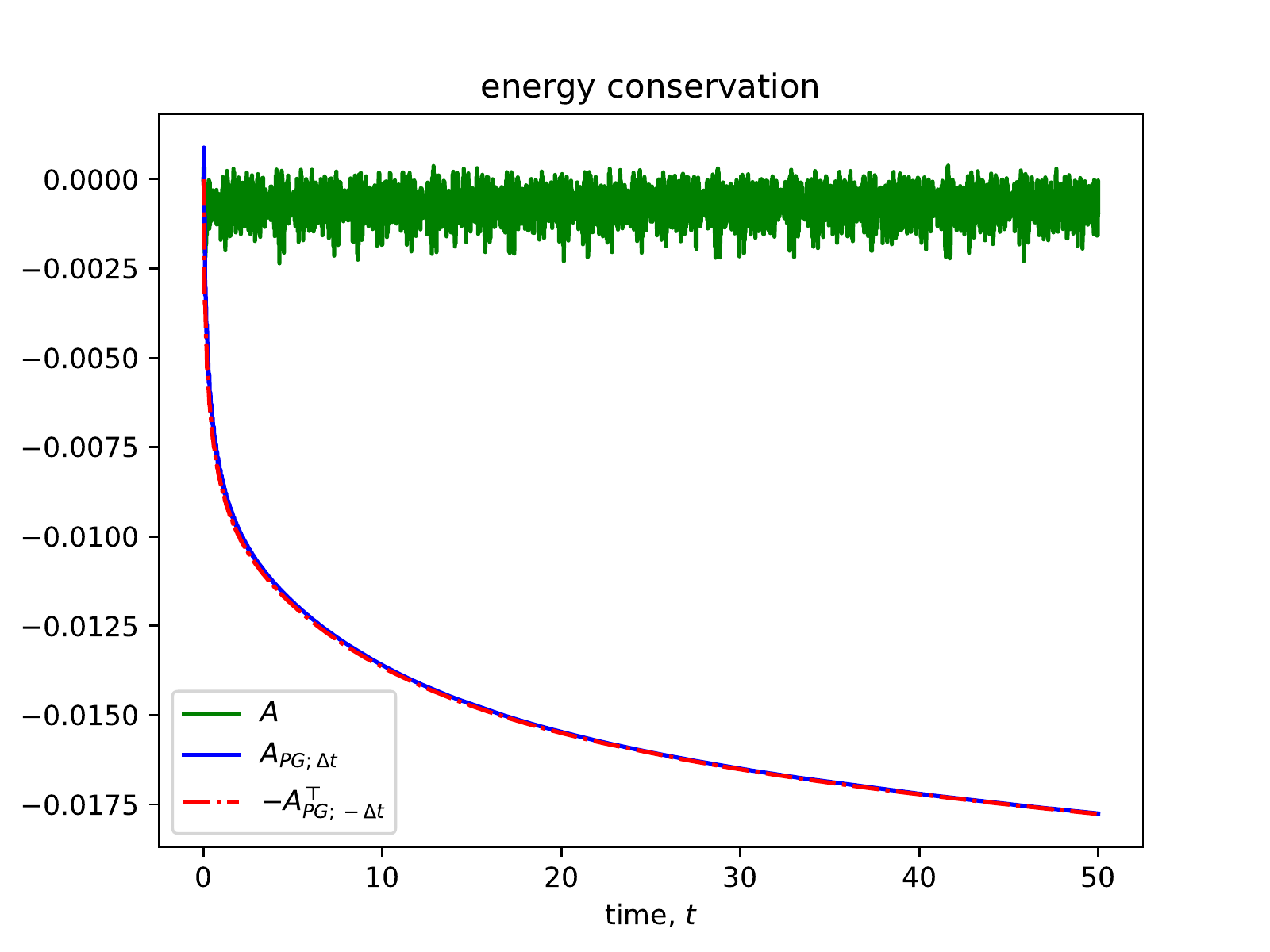}
\caption{Mass (left) and energy (right) conservation errors with time, $p=5$, 20 elements, $u=0.4$, $\Delta t=0.005$.}
\label{fig::conservation}
\end{figure}

Figure \ref{fig::advection} shows the initial and final results for the advection test. 
Note that the initial data is somewhat oscillatory since this is the projection of an analytical function with
sharp gradients onto the $\mathcal{Q}_h$ basis functions for which continuity is not enforced across element boundaries.
The upwinded test function flux form 
and downwinded trial function material form solutions are indistinguishable, however both are markedly less oscillatory than 
the original solution. The corresponding mass and energy conservation errors over 20 revolutions are given in Fig. 
\ref{fig::conservation}.
A small drift in mass conservation away from machine precision is perceptible for all three schemes, but is somewhat
greater for the upwinded schemes. Since this is present for the original advection operator, $\boldsymbol{\mathsf{A}}$ 
\eqref{eq::A}, which is known to conserve mass \cite{LPG18}, we surmise that this drift is due to the time stepping scheme.
While the energy conservation errors remain bounded for the original advection operator, the Petrov-Galerkin operators
dissipate energy.
Note that $\boldsymbol{\mathsf{A}}$ as used in \eqref{eq::centered_advection} is not itself skew-symmetric, so exact 
energy conservation is not anticipated.

\begin{figure}
\centering
\includegraphics[width=0.48\textwidth,height=0.36\textwidth]{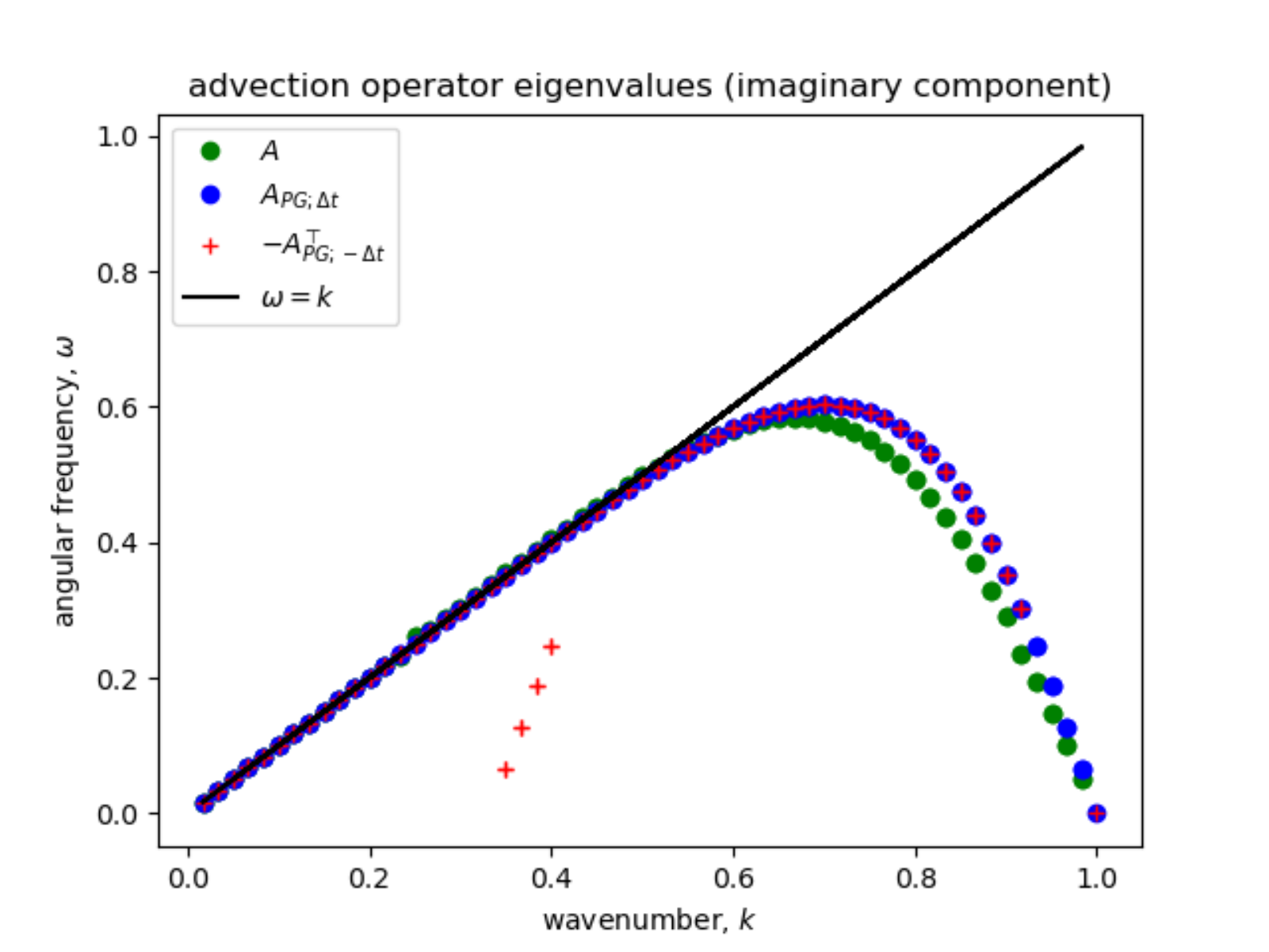}
\includegraphics[width=0.48\textwidth,height=0.36\textwidth]{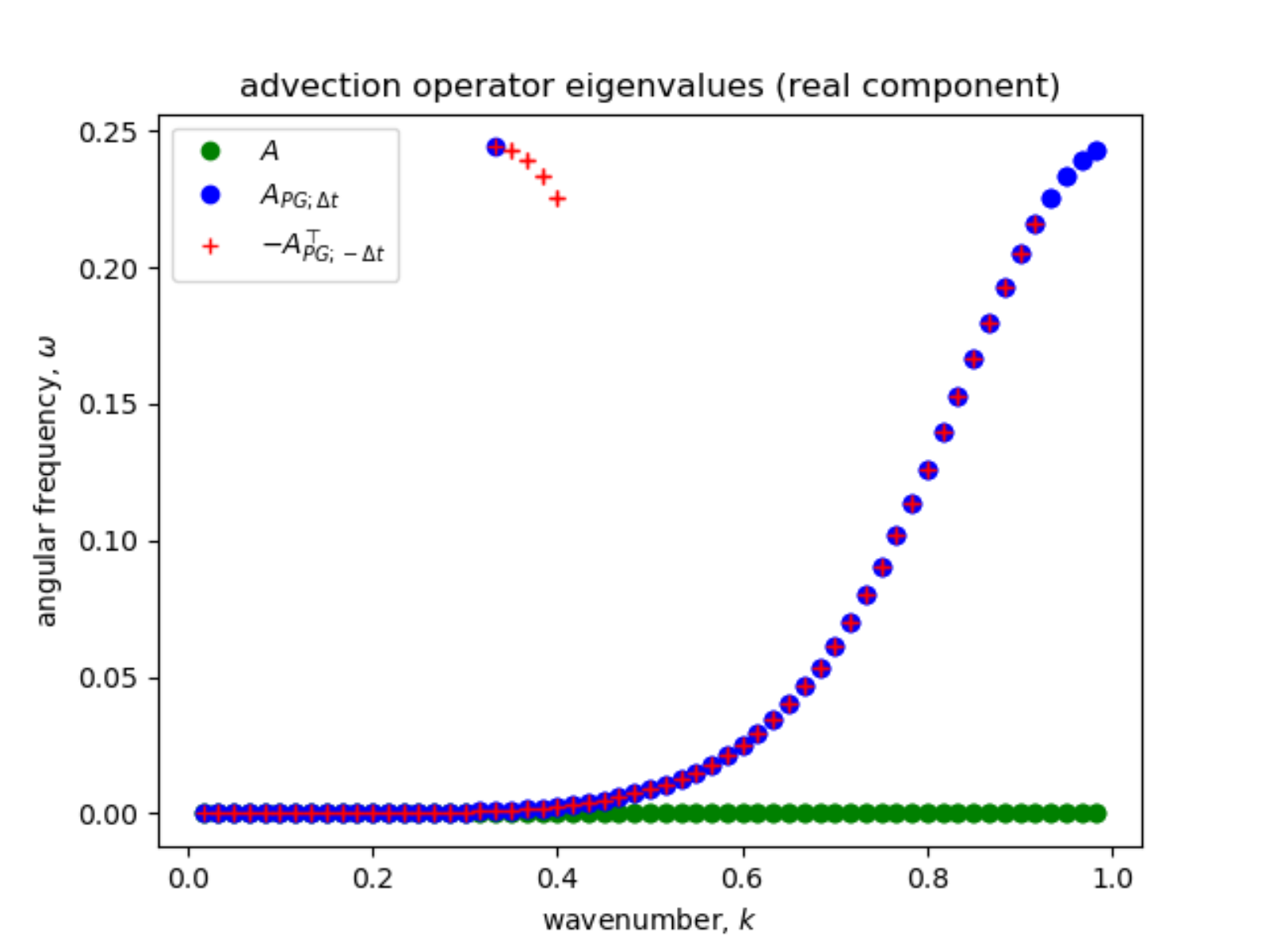}
\caption{Imaginary (left) and real (right) eigenvalues for the advection operators, 40 elements, $p=3$, $u=0.4$, $\Delta t=0.005$.}
\label{fig::eig_vals_p3}
\end{figure}

\begin{figure}
\centering
\includegraphics[width=0.48\textwidth,height=0.36\textwidth]{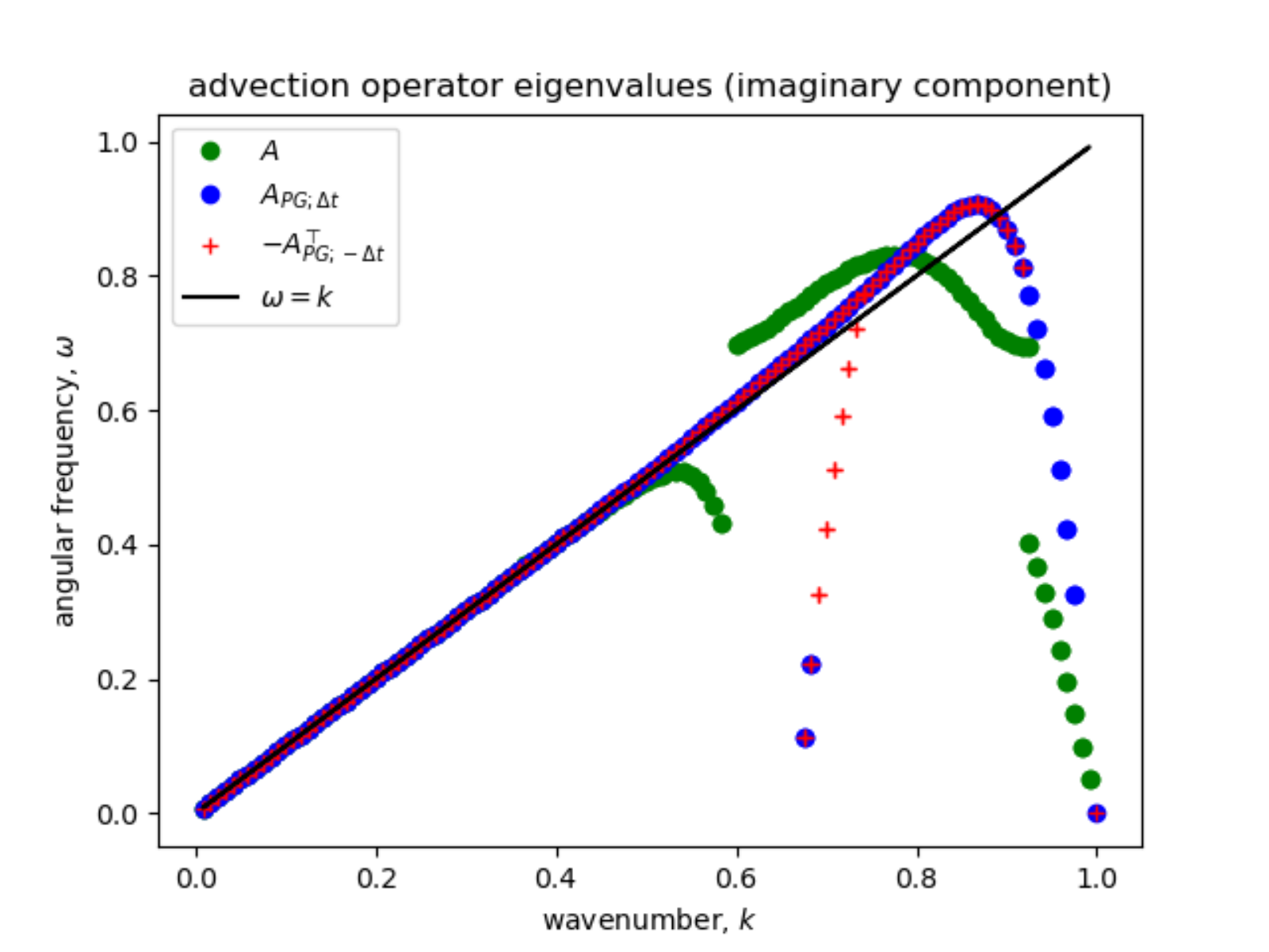}
\includegraphics[width=0.48\textwidth,height=0.36\textwidth]{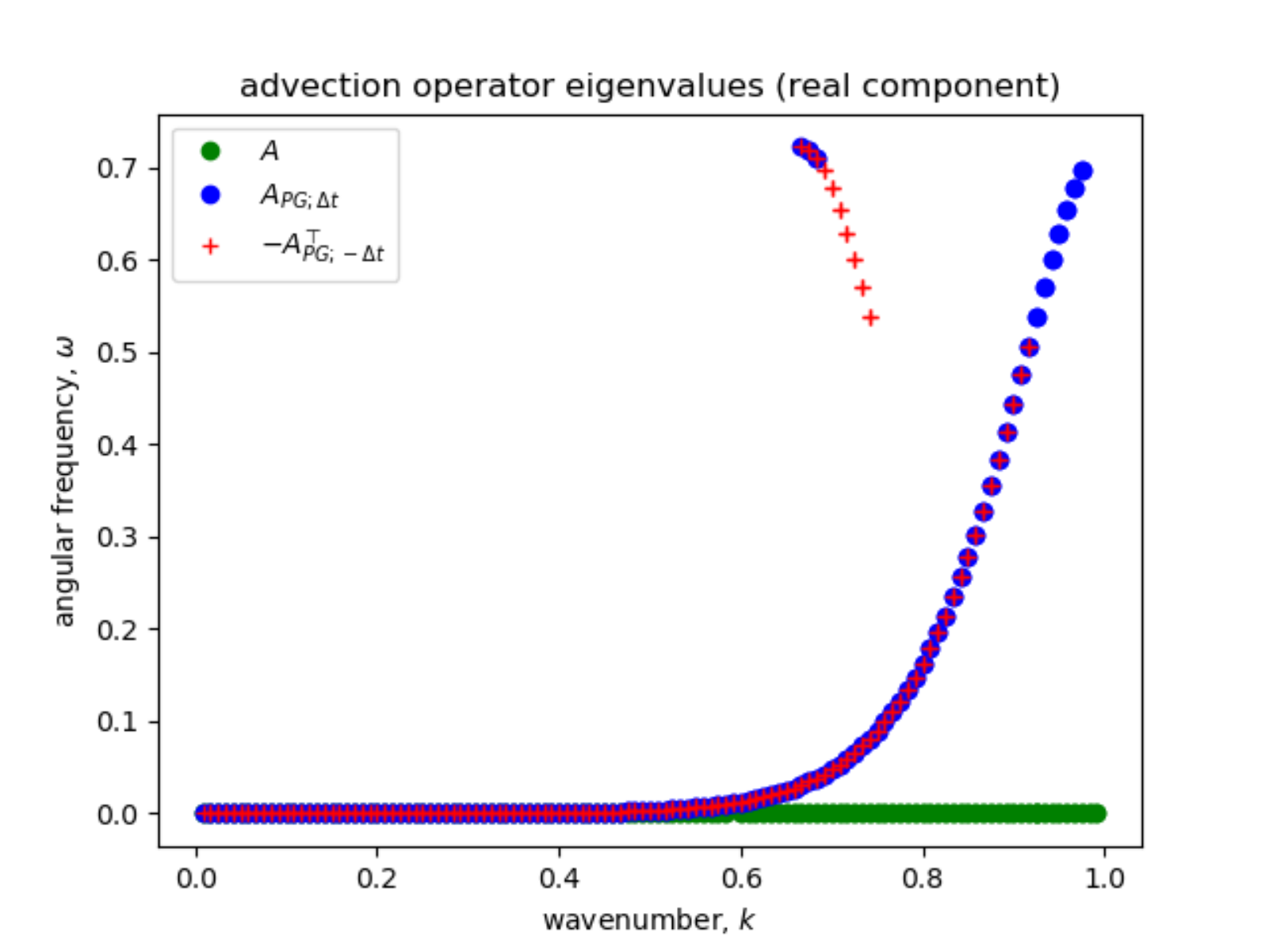}
\caption{Imaginary (left) and real (right) eigenvalues for the advection operators, 40 elements, $p=6$, $u=0.4$, $\Delta t=0.005$.}
\label{fig::eig_vals_p6}
\end{figure}

Figures \ref{fig::eig_vals_p3} and \ref{fig::eig_vals_p6} show the imaginary and real components of the dispersion relations
for the advection operators of degree $p=3$ and $p=6$ respectively with 40 elements and a time step
of $\Delta t = 0.005$. These are computed by first interpolating the eigenvectors from the $\mathcal Q_h$ space to physical
space, and then projecting these physical space solutions onto Fourier modes (for details see the appendix). 

The imaginary eigenvalues of the Petrov-Galerkin operators are closer to the analytic, dispersionless solution of $\omega=k$
(indicated by the solid black lines in Figs. \ref{fig::eig_vals_p3} and \ref{fig::eig_vals_p6}). Notably the Petrov-Galerkin
operators also seal up the spectral gaps which are characteristic of both collocated \cite{Melvin12} and mixed \cite{Eldred18} 
high order finite element discretisations. This behaviour has also been observed for hyperviscosity in the context of 
collocated spectral elements \cite{Ullrich18}. However while the real eigenvalues
are at machine precision for the original advection operator, $\boldsymbol{\mathsf{A}}$, indicating purely hyperbolic 
advection, these are non-zero for the 
Petrov-Galerkin formulations, indicating dissipative solutions biased towards higher wave numbers. The dissipation profile
of these operators steepens with polynomial degree, in analogy to a higher power viscosity operator.

\begin{figure}
\centering
\includegraphics[width=0.48\textwidth,height=0.36\textwidth]{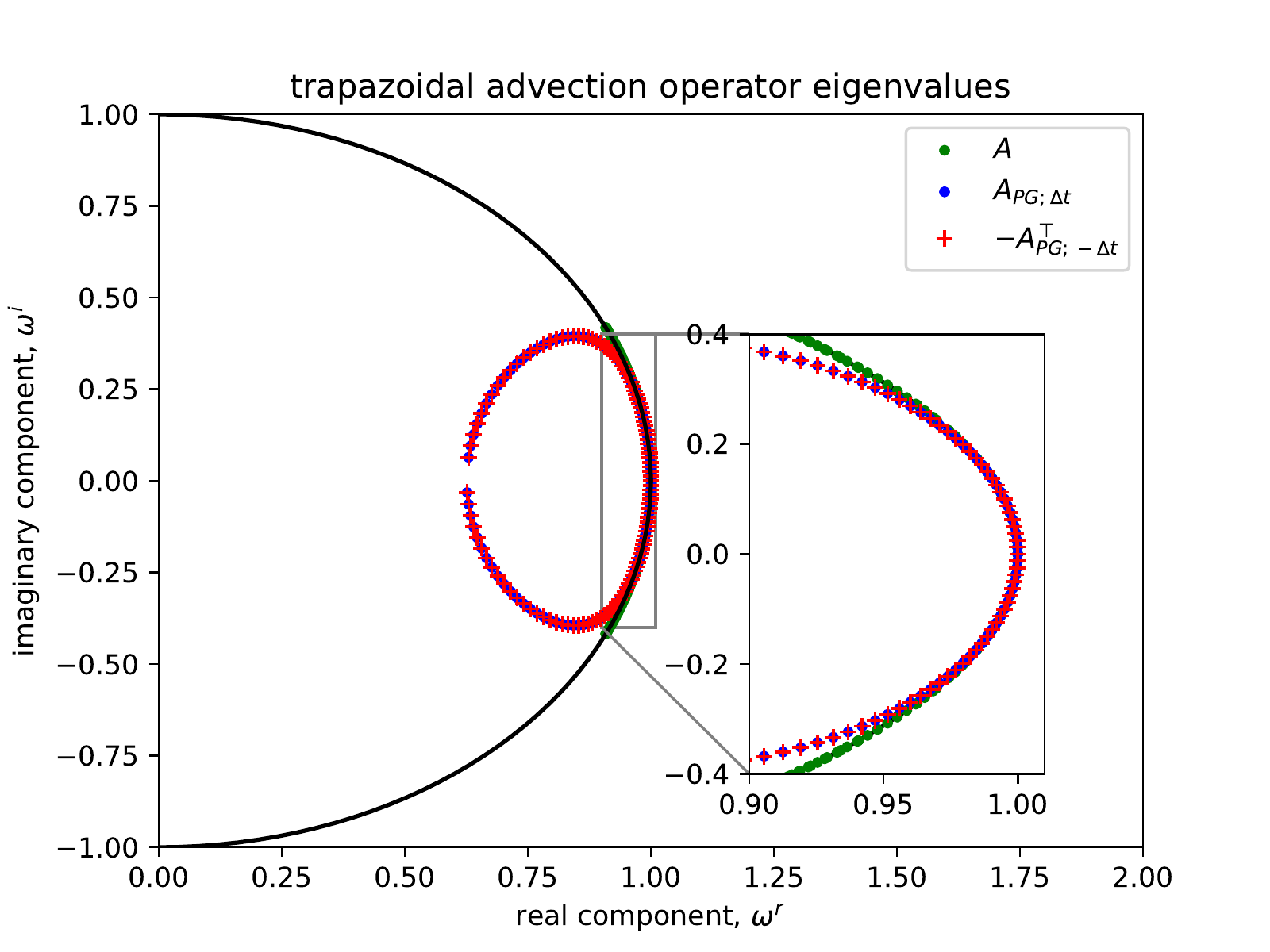}
\includegraphics[width=0.48\textwidth,height=0.36\textwidth]{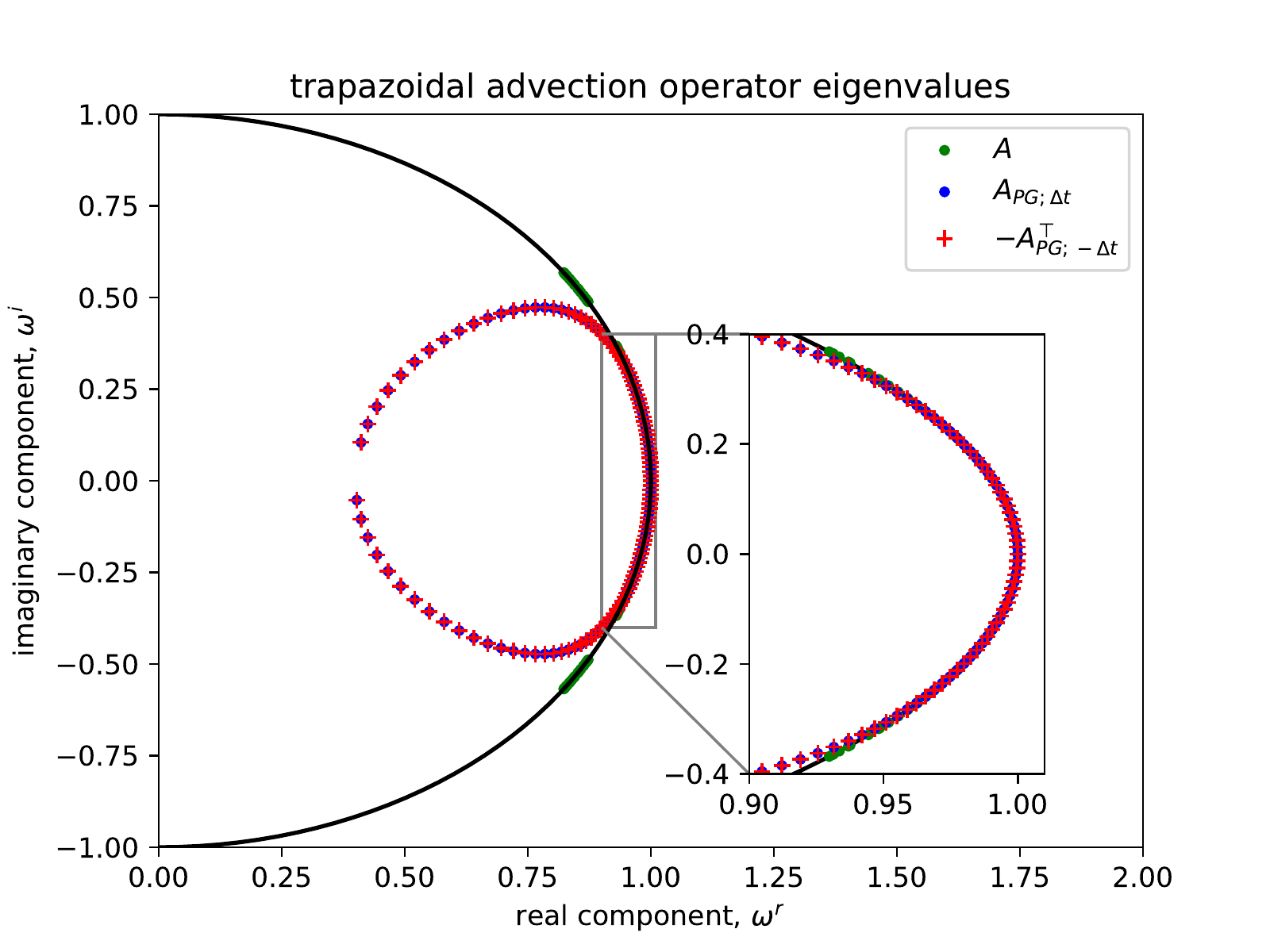}
\caption{Time centered advection operator eigenvalues (real), $p=3$, 40 elements (left), and 
$p=6$, 20 elements (right), $u=0.4$, $\Delta t=0.005$. The black line indicates the unit circle.}
\label{fig::stab_region}
\end{figure}

Figure \ref{fig::stab_region} shows the real component of the eigenvalues, $\omega^r$ against the imaginary, 
$\omega^i$ for the centered time stepping operator 
$(\boldsymbol{\mathsf M} + 0.5\Delta t\boldsymbol{\mathsf A})^{-1}(\boldsymbol{\mathsf M} - 0.5\Delta t\boldsymbol{\mathsf A})$
(and similarly for $\boldsymbol{\mathsf A}_{PG;\Delta t}$ and $-\boldsymbol{\mathsf A}_{PG;-\Delta t}^{\top}$), as given 
in \eqref{eq::centered_advection}, for a time step of $\Delta t = 0.005$. The eigenvalues determined by evaluating the time
stepping operator with respect to $\boldsymbol{\mathsf A}$ sit on the 
unit circle, indicating the neutral stability of this formulation. This is also reflected in the bounded energy conservation
errors, as shown in Fig. \ref{fig::conservation}. The dissipative nature of the Petrov-Galerkin operators is reflected by the 
fact that the eigenvalues for the time integration operator sit \emph{inside} the unit circle.

\begin{figure}
\centering
\includegraphics[width=0.48\textwidth,height=0.36\textwidth]{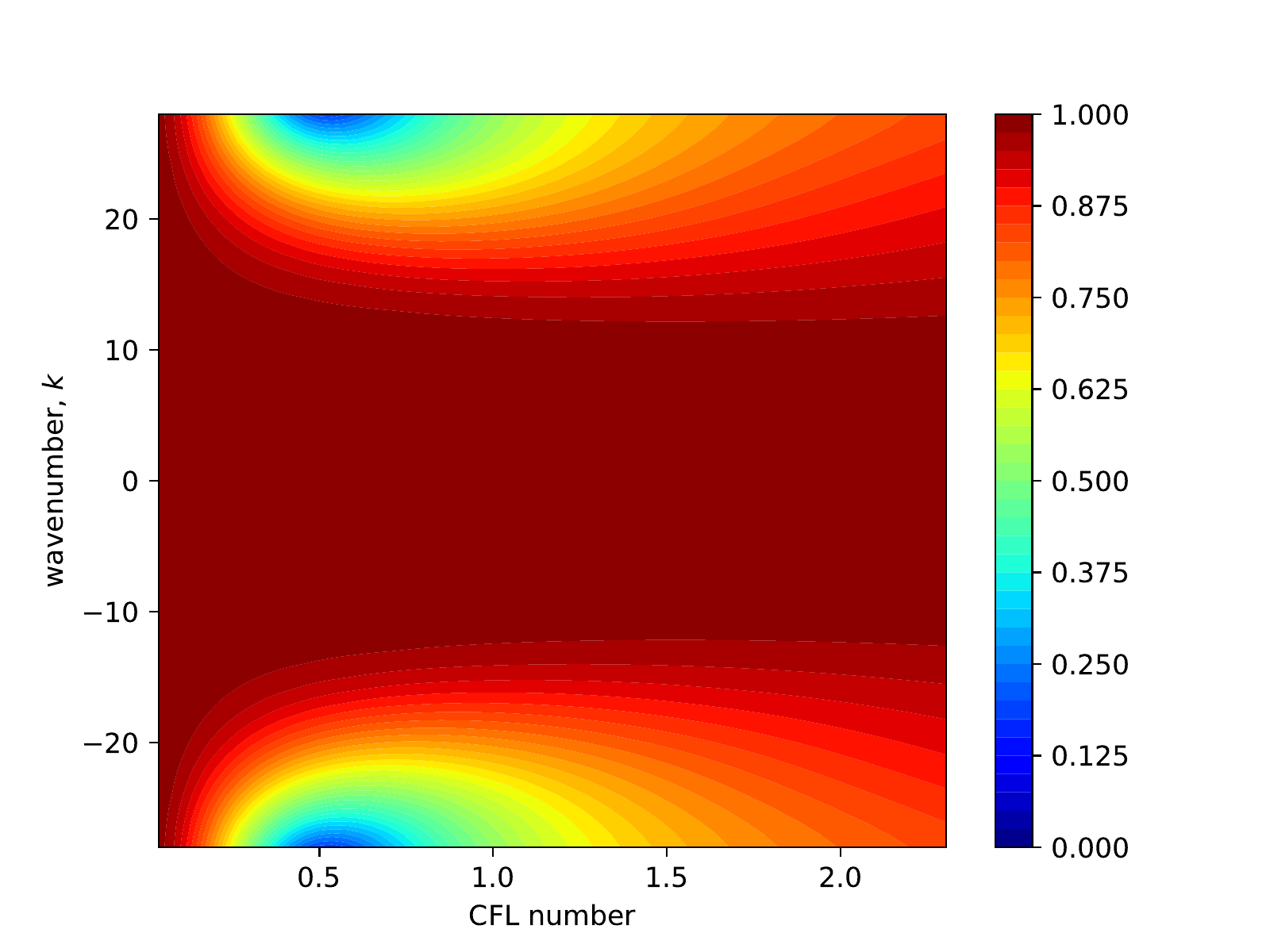}
\includegraphics[width=0.48\textwidth,height=0.36\textwidth]{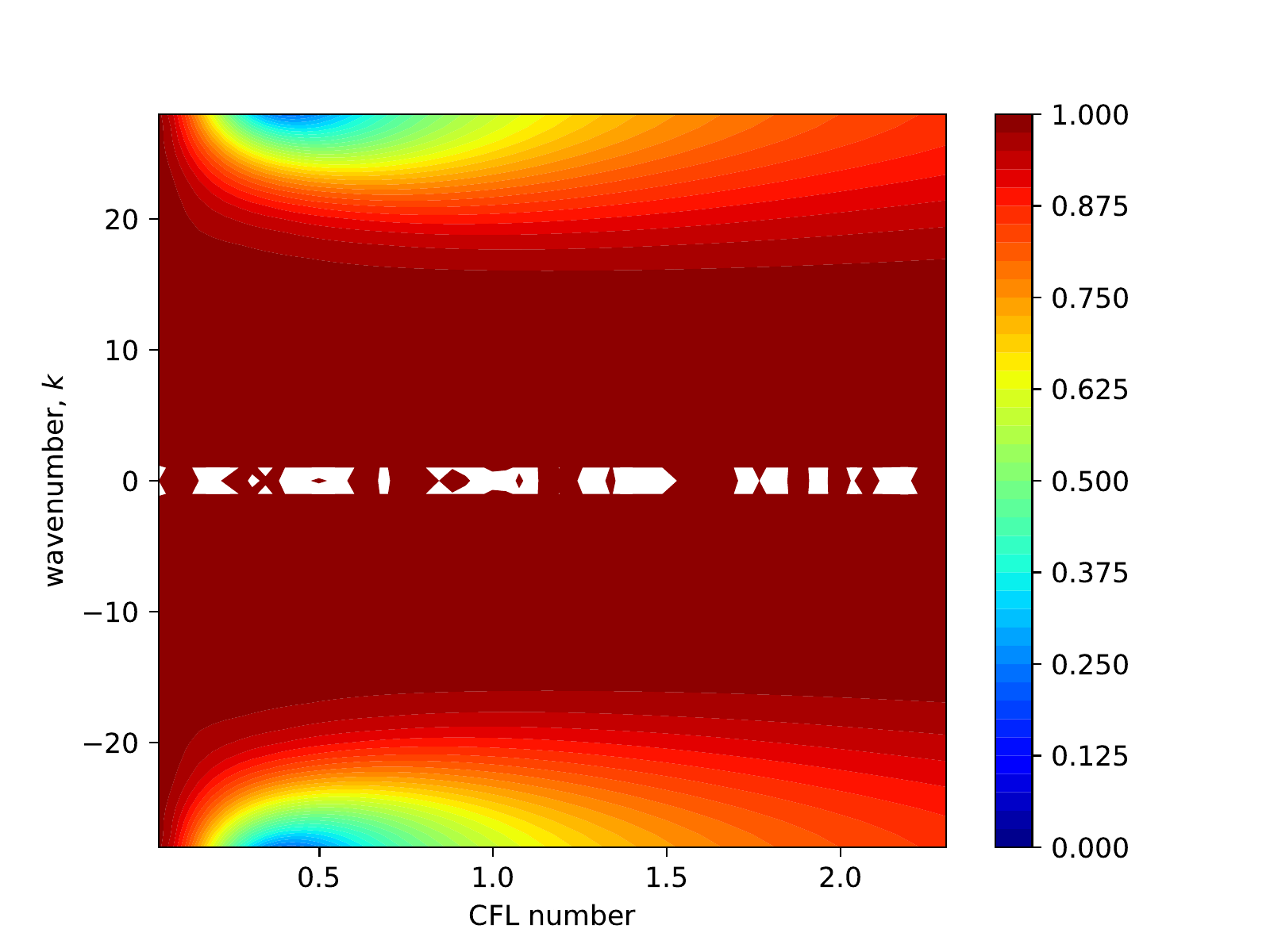}
\caption{Absolute eigenvalue of the time centered flux form Petrov-Galerkin advection operator, $|\omega|$, as a function of CFL number, 
$\Delta t|\boldsymbol{u}|n_ep/L$ and Fourier wavenumber, $k$; $p=3$, $n_e=20$ (left) and $p=6$, $n_e=10$ (right).}
\label{fig::A_stability}
\end{figure}

The A-stability of the upwinded flux form time stepping operator,
$(\boldsymbol{\mathsf M} + 0.5\Delta t\boldsymbol{\mathsf{A}}_{PG;\Delta t})^{-1}(\boldsymbol{\mathsf M} - 0.5\Delta t\boldsymbol{\mathsf{A}}_{PG;\Delta t})$,
is determined from the magnitude of its eigenvalues, $|\omega|$. These are plotted against an approximate CFL number,
$\Delta t|\boldsymbol{u}|n_ep/L$ (which does not account for the smaller distance between the GLL nodes near the 
element boundaries), and Fourier wavenumber, $k$ for polynomials of degree $p=3$ and $p=6$ in Fig. 
\ref{fig::A_stability}. While dissipation increases for larger absolute wavenumbers $|k|$, as indicated by smaller values of $|\omega|$, 
it is greatest for moderate CFL numbers. Since the time centered advection operator is unconditionally
stable, nowhere is $|\omega| > 1$. Additionally, for larger CFL numbers the upwinded operators fail to 
preserve their theoretical rate of convergence.

\begin{figure}
\centering
\includegraphics[width=0.48\textwidth,height=0.36\textwidth]{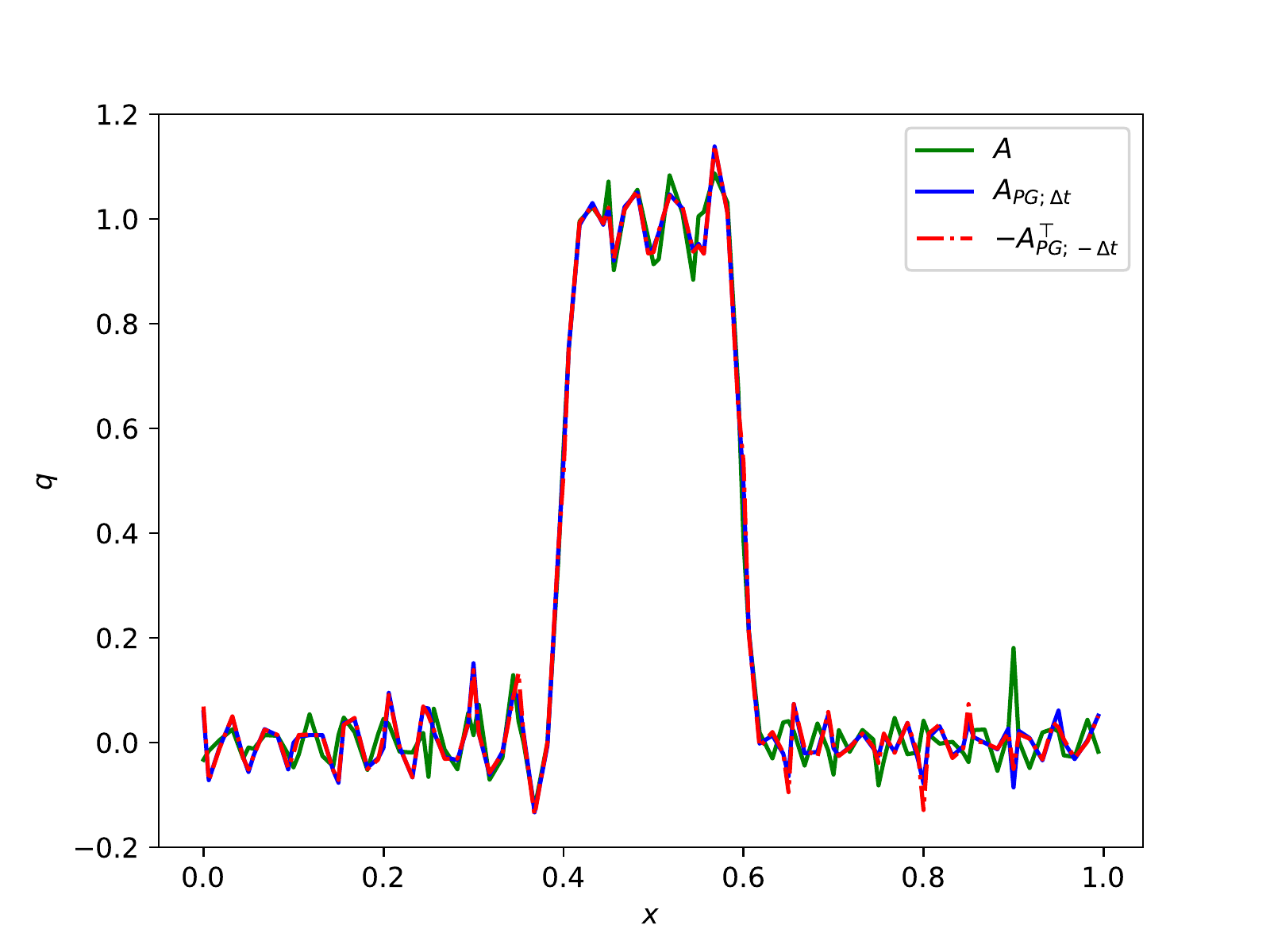}
\includegraphics[width=0.48\textwidth,height=0.36\textwidth]{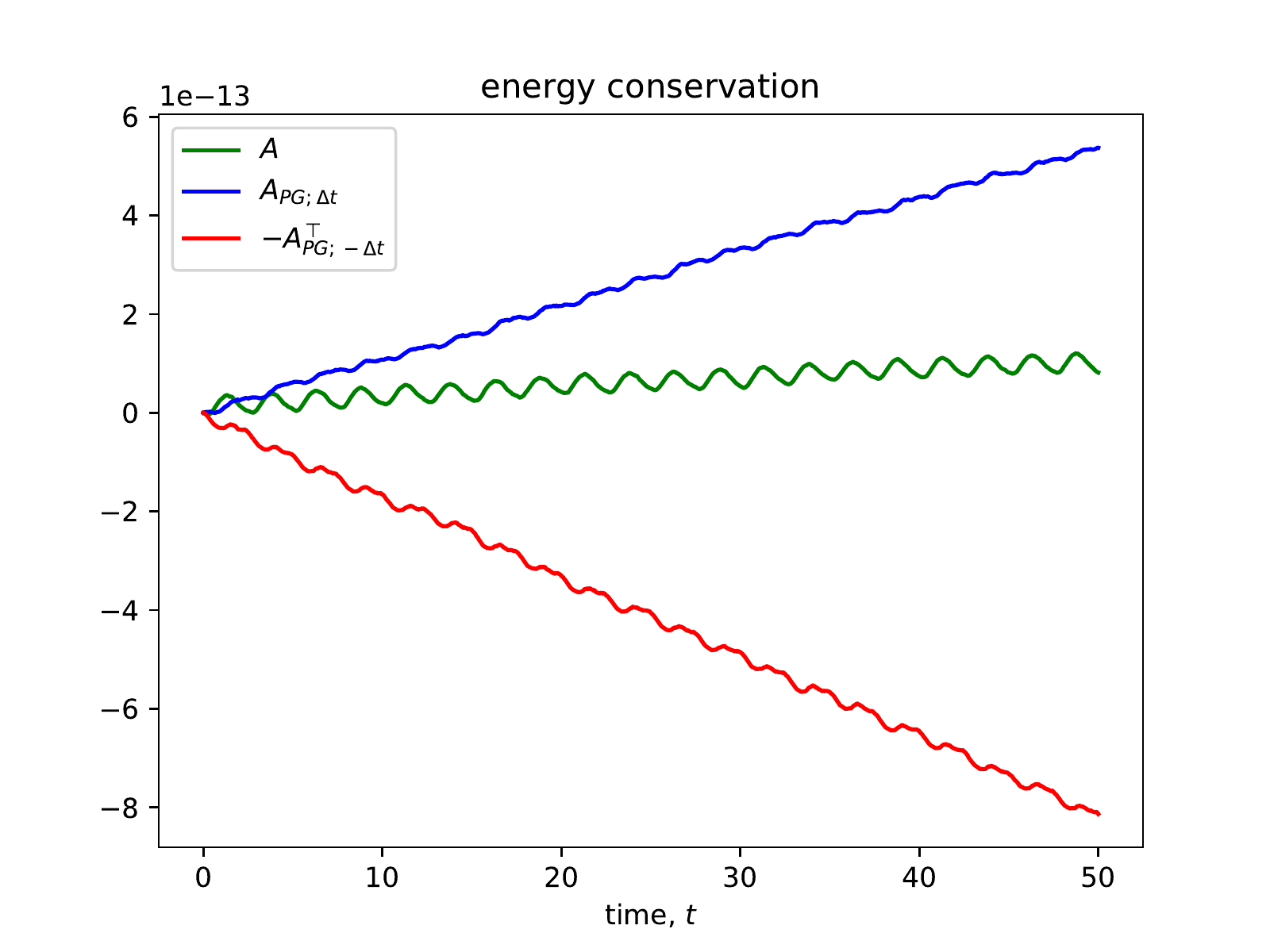}
\caption{Skew-symmetric formulation: final state after one revolution (left) and energy 
conservation errors (right), $p=5$, 20 elements, $u=0.4$, $\Delta t=0.005$.}
\label{fig::skew_sym}
\end{figure}

As for the original advection operator, $\boldsymbol{\mathsf A}$, both forms of the Petrov Galerkin 
advection operator may be used to construct skew-symmetric formulations, as given in \eqref{adv_ss}. In each 
case there is a small drift away from machine precision in the energy conservation as shown in Fig. 
\ref{fig::skew_sym}, which is also observed for the mass conservation in Fig. \ref{fig::conservation} and is
due to the mass conservation errors in the time stepping scheme. These skew-symmetric formulations  directly 
result from the cancellation of the upwinding
contributions of the Petrov-Galerkin operators such that the solutions are once again oscillatory and there
is no discernible benefit to the use of these formulations over the original advection operator.

\section{Extension to multiple dimensions: Advection on the sphere}

Having validated the new advection scheme for the one dimensional case, the scheme is extended
to the case of passive advection on the surface of the sphere. The compatibility property between 
the discrete spaces $U_h\subset H(\mathrm{div},\Omega^2)$ and $Q_h\subset L^2(\Omega^2)$ \eqref{H_div_L_2} 
is also satisfied in the two dimensional domain, $\Omega^2 \subset \mathbb{R}^2$, as is the strong form mapping of the 
incidence matrix, $\boldsymbol{\mathsf{E}}^{2,1}$
between these spaces and the weak form adjoint relationship to the discrete gradient operator \eqref{disc_grad}. 
However in multiple dimensions the bases that span these function spaces are constructed from tensor product 
combinations of the nodal and edge polynomials \eqref{eq::nodal_edge}. These are given respectively 
for the two dimensional basis functions $\boldsymbol{\beta}_h\in U_h$ and $\gamma_h\in Q_h$ for
polynomial degree $p$ as
\begin{subequations}\label{beta_gamma}
\begin{align}
\boldsymbol{\beta}_k(\xi,\eta) &:= 
\begin{cases}
l_i(\xi)e_j(\eta)\boldsymbol{e}_{\xi}\quad &
\text{if $k$ even,$\quad$ with $i = 0,\dots,p\qquad j = 0,\dots,p-1\quad k=2(j(p+1) + i)$}\\
e_i(\xi)l_j(\eta)\boldsymbol{e}_{\eta}\quad &
\text{if $k$ odd,$\quad$ with $i = 0,\dots,p-1\qquad j = 0,\dots,p\quad k=2(jp+i)+1$},
\end{cases}\label{basis_beta}
\\
\gamma_k(\xi,\eta) &:= e_i(\xi)e_j(\eta),\quad i,j = 0,\dots,p-1\quad k = jp + i.
\end{align}
\end{subequations}
for which
\begin{equation}
U_h = \mathrm{span}\{\boldsymbol{\beta}_0(\xi,\eta),\dots,\boldsymbol{\beta}_{2p(p+1)-1}(\xi,\eta)\},\quad
Q_h = \mathrm{span}\{\gamma_0(\xi,\eta),\dots,\gamma_{p\times p-1}(\xi,\eta)\},
\end{equation}
where $(\xi,\eta)$ are the local coordinates in the dimensions $\boldsymbol{e}_{\xi}$, 
$\boldsymbol{e}_{\eta}$ respectively within the two dimensional canonical domain $[-1,1]\times[-1,1]$.

In order to upwind the $U_h$ bases, these are evaluated at downstream locations according to
\eqref{xi_downwind} only in dimensions for which $C^0$ continuity is enforced, such that the 
upwinded form of \eqref{basis_beta} is given as
\begin{equation}
\boldsymbol{\beta}^u_k(\xi,\eta) := 
\begin{cases}
l_i(\xi^d)e_j(\eta)\boldsymbol{e}_{\xi}\quad &
\text{if $k$ even,$\quad$ with $i = 0,\dots,p\qquad j = 0,\dots,p-1\quad k=2(j(p+1) + i)$}\\
e_i(\xi)l_j(\eta^d)\boldsymbol{e}_{\eta}\quad &
\text{if $k$ odd,$\quad$ with $i = 0,\dots,p-1\qquad j = 0,\dots,p\quad k=2(jp+i)+1$}.
\end{cases}
\end{equation}
The two dimensional upwinded flux form advection equation is then given as
\begin{equation}\label{adv_flux_2d}
\langle\gamma_i,\gamma_j\rangle_{\Omega^2}\dot{\hat{q}}_j + 
\langle\gamma _i,\gamma_k\rangle_{\Omega^2}\mathsf{E}_{k,l}^{2,1}\langle\boldsymbol{\beta}^u_m,\boldsymbol{\beta}_l\rangle^{-1}_{\Omega^2}
\langle\boldsymbol{\beta}^u_m\cdot\boldsymbol{u}_h,\gamma_n\rangle_{\Omega^2}\hat{q}_n,\qquad\forall\gamma_i\in Q_h.
\end{equation}

In multiple dimensions on parallel machines it is not practical to explicitly construct a mass
matrix inverse for the $U_h$ space. Consequently we use an explicit third order Runge-Kutta scheme
in place of the centered time integration scheme used in the previous section. This integrator 
takes the form \cite{Durran10}
\begin{subequations}\label{rk_3}
\begin{align}
q_h^{(1)} &= q_h^n - \Delta ty(\boldsymbol{u}_h^n,q_h^n),\\
q_h^{(2)} &= q_h^n - \frac{\Delta t}{4}\Big(y(\boldsymbol{u}_h^n,q_h^n) + y(\boldsymbol{u}_h^{n+1},q_h^{(1)})\Big),\\
q_h^{n+1} &= q_h^n - \frac{\Delta t}{6}\Big(y(\boldsymbol{u}_h^n,q_h^n) + y(\boldsymbol{u}_h^{n+1},q_h^{(1)}) + 4y(\boldsymbol{u}_h^{n+1/2},q_h^{(2)})\Big),
\end{align}
\end{subequations}
where $y(\boldsymbol{u}_h,q_h)$ represents the second term in \eqref{adv_flux_2d}.

\begin{figure}
\centering
\includegraphics[width=0.48\textwidth,height=0.36\textwidth]{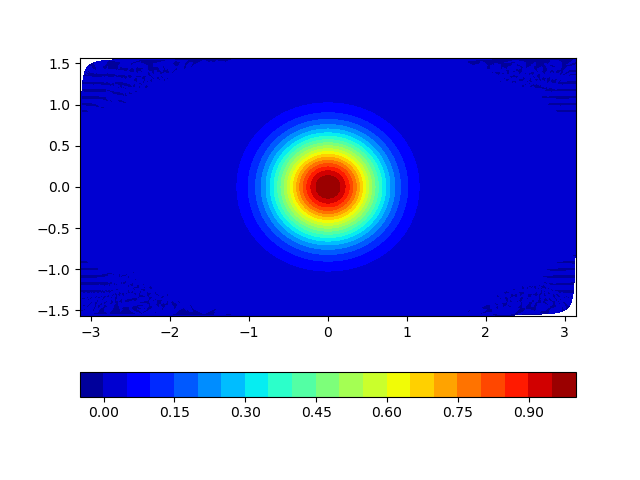}
\includegraphics[width=0.48\textwidth,height=0.36\textwidth]{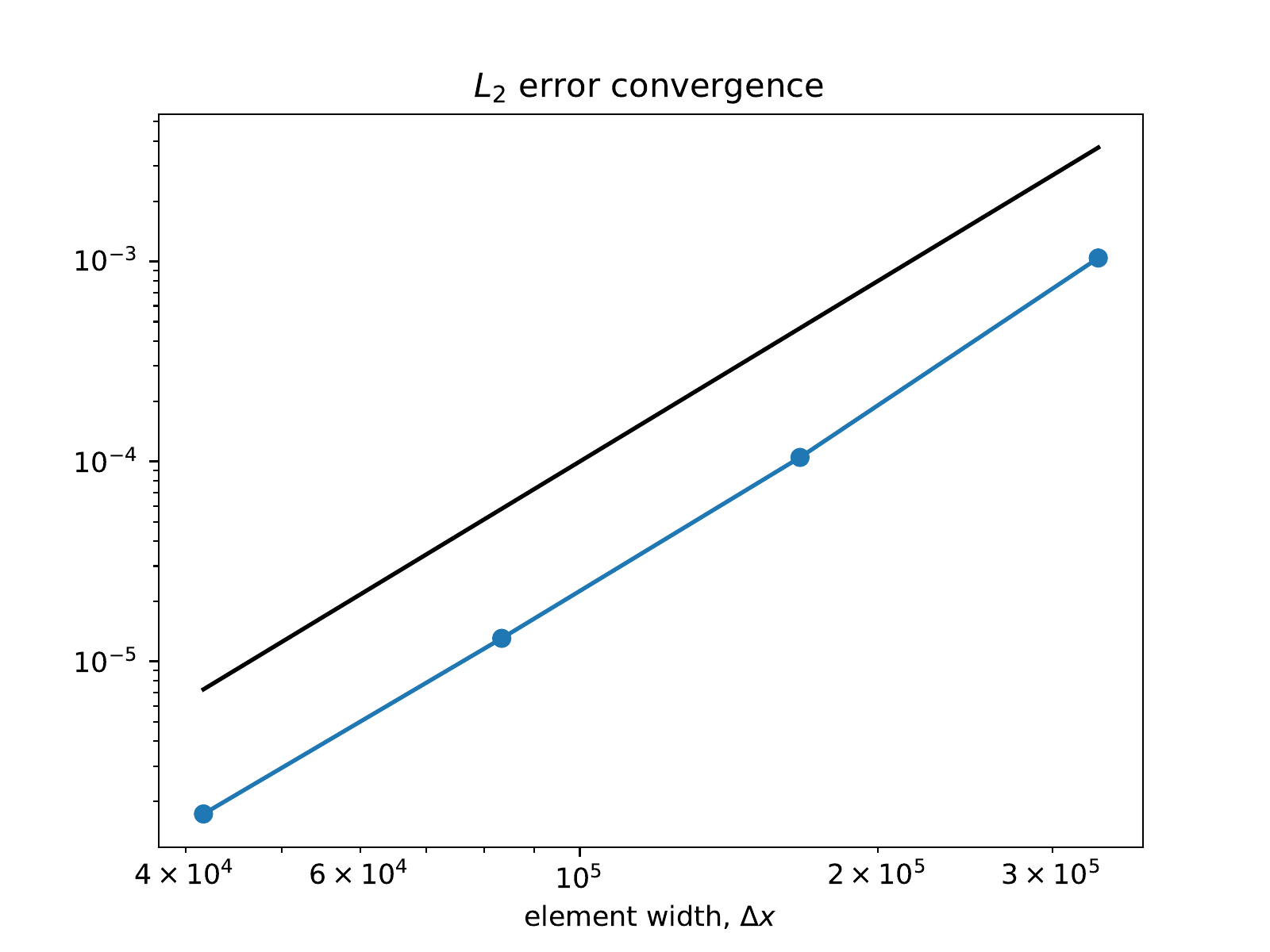}
\caption{Solid body advection on the sphere: tracer field at final time (left) and $L_2$ error convergence (right).
Solid line indicates a gradient of $\Delta x^3$.}
\label{fig::2d_advection_0}
\end{figure}

The advection equation as described above is implemented using a cubed sphere discretisation with
a physical radius of 6371220.0m, the details of which can be found in \cite{LP18}. The first test
on the sphere is one of solid body rotation, with a prescribed velocity field of
$\boldsymbol{u}(\theta,\phi) = (20\cos(\phi)\boldsymbol{e}_{\theta},0\boldsymbol{e}_{\phi})$,
where $\boldsymbol{e}_{\theta}$ and $\boldsymbol{e}_{\phi}$ are the unit vectors in the zonal and
meridional directions respectively, and an initial Gaussian tracer field of the form 
$q(t=0) = \cos(\phi)e^{-9\theta^2 - 225\phi^2}$. The test is run for a single revolution using 
elements of degree $p=3$ with $384\times n_e/10$ time steps, where $n_e$ is the number of elements in
each dimension on each of the six panels of the cubed sphere, for which $n_e = 10,20,40,80$.
Figure \ref{fig::2d_advection_0} shows the final tracer distribution for $n_e=20$, as well as the 
convergence of the normalised $L^2$ error of the tracer field after one revolution. Using both a
third order time integration scheme \eqref{rk_3} and elements of degree $p=3$ the errors converge
at third order.

The second test on the sphere is for the advection of two cosine bells within a time varying velocity field
that involves both zonal transport and shear deformation, with the tracer field returning to its original
position over a period of $T=12$ days. The specific details of the tracer and velocity field
configurations are described in \cite{Lauritzen12}. The test is run with $40\times 40$ elements 
of degree $p=3$ on each face of the cubed sphere, with a time step of $\Delta t = 129.6\mathrm{s}$.

\begin{figure}
\centering
\includegraphics[width=0.48\textwidth,height=0.36\textwidth]{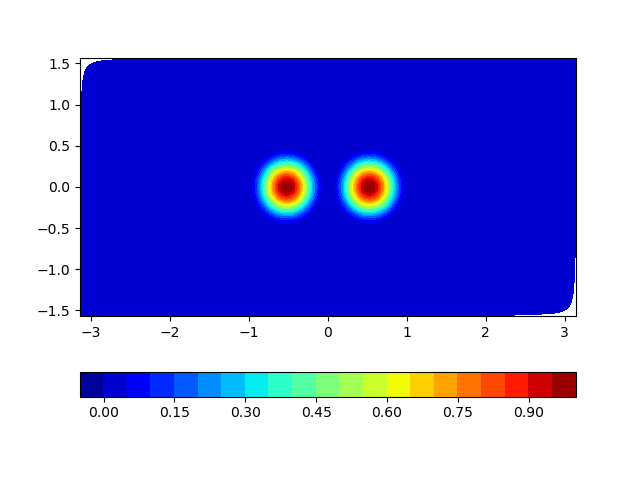}
\includegraphics[width=0.48\textwidth,height=0.36\textwidth]{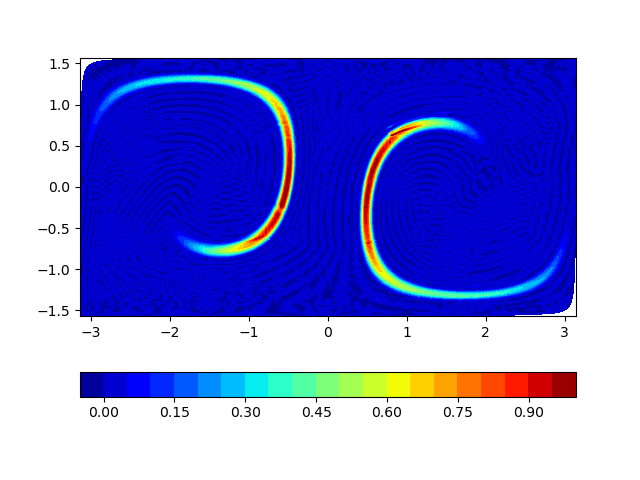}
\caption{Deformational flow on the sphere: tracer field at days 0 (left) and 6 (right).}
\label{fig::2d_advection_1}
\end{figure}

\begin{figure}
\centering
\includegraphics[width=0.48\textwidth,height=0.36\textwidth]{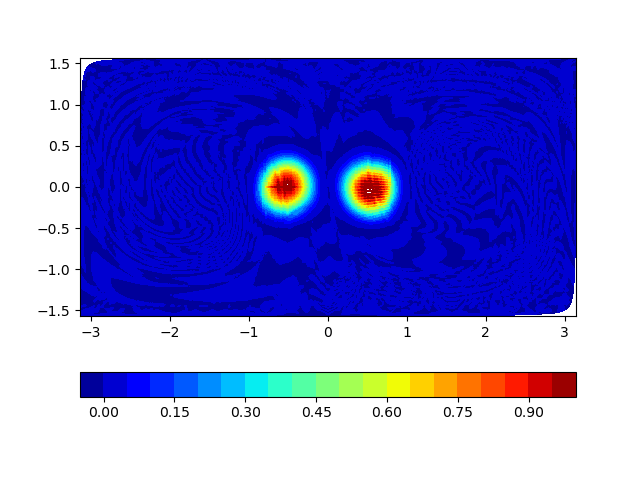}
\includegraphics[width=0.48\textwidth,height=0.36\textwidth]{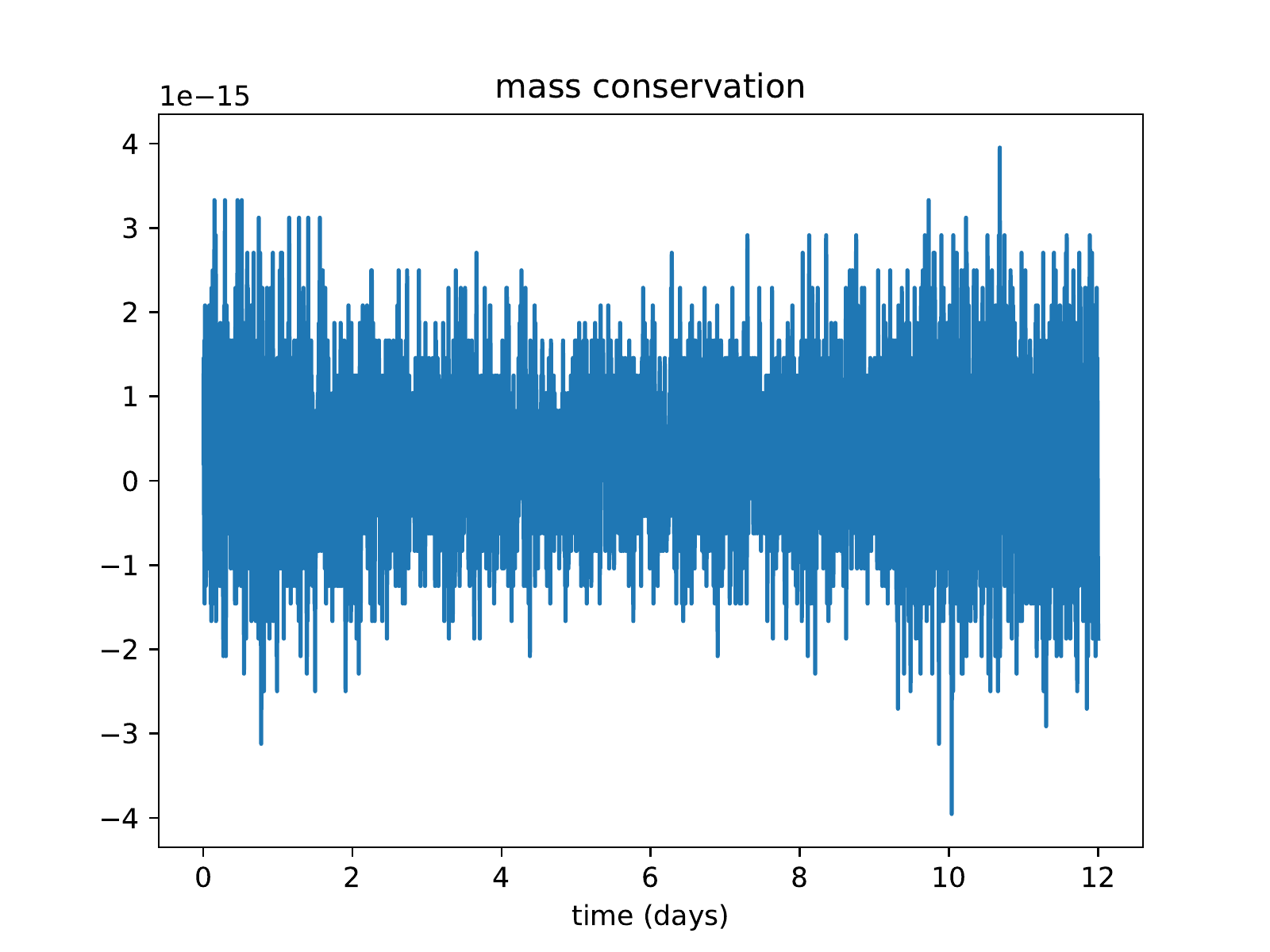}
\caption{Deformational flow on the sphere: tracer field at day 12 (left) and mass conservation error (right).}
\label{fig::2d_advection_2}
\end{figure}

Figures \ref{fig::2d_advection_1} and \ref{fig::2d_advection_2} show the evolution of the tracer field at
days 0, 6 and 12, as well as the mass conservation errors. While there is a small drift in the mass conservation
errors for the one dimensional test case owing to the semi-implicit time integration scheme, as observed in 
Fig. \ref{fig::conservation}, no drift is evident in the two dimensional advection problem, for which an explicit
integration scheme \eqref{rk_3} is used. Notably, in the absence of upwinding the tracer field becomes completely
incoherent over the course of the simulation, such that the upwinded correction is essential for this test case. 
The solution errors do not converge at the expected rate for this 
test case, most likely on account of the time varying velocity field and the first order integration of the characteristics
used to compute the test function nodal departure points. In order to compute these at higher order the velocity 
field would need to be interpolated at intermediate locations. Currently the parallel communication of data within
the code is only supported along element boundaries, and not within the internal element locations needed to allow for this.

\section{Application to potential vorticity stabilisation in shallow water on the sphere}

In order to demonstrate the usefulness of the Petrov Galerkin upwinding described 
above this is applied to the diagnosis and interpolation of the potential vorticity for the rotating 
shallow water equations on the sphere. 
In additional to the discrete bases spanning the $U_h$ and $Q_h$ spaces described in the previous
section \eqref{beta_gamma}, an additional space, $W_h\subset H(\mathrm{rot},\Omega^2)$ is introduced,
for which the basis functions are given as

\begin{equation}
\alpha_k(\xi,\eta) := l_i(\xi)l_j(\eta),\quad i,j=0\dots p\quad k = j(p+1)+i.
\end{equation}
As well as the \emph{div} compatibilty property \eqref{H_div_L_2}, there is also a compatible 
mapping between the discrete spaces $W_h$ and $U_h$ via the \emph{rot} operator. Together these give 
the sequence

\begin{equation}\label{H_rot_H_div_L_2}
\mathbb{R}\longrightarrow W_h \stackrel{\nabla^{\perp}}{\longrightarrow} U_h \stackrel{\nabla\cdot}{\longrightarrow} Q_h \longrightarrow 0.
\end{equation}

Since the bases for the $W_h$ space are $C^0$ continuous in both dimensions, the downwinded variant 
is given for $(\xi^u,\eta^u)$ as defined in \eqref{xi_upwind} as

\begin{equation}
\alpha_k^d(\xi,\eta) := l_i(\xi^u)l_j(\eta^u),\quad i,j=0\dots p\quad k = j(p+1)+i.
\end{equation}

In addition to the mixed mimetic spectral element spatial discretisation \cite{LP18}, the model uses a
semi-implicit time integration scheme \cite{WCB20}, that together allow for the exact balance of 
energy exchanges, and the exact conservation of mass, vorticity and energy in space and time.
No dissipation is applied to the model, except that in the second simulation the upwinded 
stabilisation is applied by sampling the potential vorticity field trial functions at downstream locations 
in a two dimensional analogue of \eqref{B_PG}. Just as there is an adjoint relation between the discrete
\emph{div} and \emph{grad} operators \eqref{disc_grad}, a similar adjoint relationship exists between 
the discrete \emph{rot} and \emph{curl} operators \cite{Kreeft13,LPG18}. The downwinded potential vorticity field 
is therefore diagnosed in the discrete weak form analogue of 
$q = (\nabla\times\boldsymbol u + f)/h$, where $q$ is the potential vorticity, $f$ is the Coriolis term and
$h$ is the fluid depth, as
\begin{equation}\label{vort_eqn_up}
\langle\alpha_i,h_h\alpha_j^{d}\rangle_{\Omega^2}\hat{q}_j^{PG} = 
-(\mathsf{E}_{k,j}^{1,0})^{\top}\langle\boldsymbol{\beta}_k,\boldsymbol{\beta}_l\rangle_{\Omega^2}\hat{u}_l
+ \langle\alpha_i,\alpha_j\rangle_{\Omega^2}\hat{f}_j,
\qquad\forall\alpha_i\in {W}_h,
\end{equation}
where $\boldsymbol{\mathsf{E}}^{1,0}$ is the incidence matrix representing the discrete strong form \emph{rot} operator \cite{LPG18,LP18}. 
The upwinded potential vorticity is then coupled to the shallow water system through the vector invariant 
form of the momentum equation 
$\dot{\boldsymbol u} = -q\times\boldsymbol{M} - \nabla\Phi$, where $\boldsymbol{M} = h{\boldsymbol{u}}$
is the mass flux and $\Phi = |u|^2/2 + gh$ is the Bernoulli function and $\dot{h} = -\nabla\cdot\boldsymbol{M}$,
\cite{LP18} for all $\boldsymbol{\beta}_i\in {U}_h$, $\gamma_i\in Q_h$ as
\begin{equation}\label{sw_mom_eq}
\begin{bmatrix}
\langle\boldsymbol{\beta}_i,\boldsymbol{\beta}_j\rangle_{\Omega^2}\hat{u}_j^{n+1} \\
\langle\gamma_i,\gamma_j\rangle_{\Omega^2}\hat{h}_j^{n+1} 
\end{bmatrix} = 
\begin{bmatrix}
\langle\boldsymbol{\beta}_i,\boldsymbol{\beta}_j\rangle_{\Omega^2}\hat{u}_j^{n} \\
\langle\gamma_i,\gamma_j\rangle_{\Omega^2}\hat{h}_j^{n} 
\end{bmatrix} - \Delta t
\begin{bmatrix}
\langle\boldsymbol{\beta}_i,\overline{q}_h^{d}\times\boldsymbol{\beta}_k\rangle_{\Omega^2} &
-(\mathsf{E}^{2,1}_{l,i})^{\top}\langle\gamma_l,\gamma_m\rangle_{\Omega^2} \\
\langle\gamma_i,\gamma_j\rangle_{\Omega^2}\mathsf{E}^{2,1}_{j,k} & \mathsf{0}
\end{bmatrix}
\begin{bmatrix}
\overline{\hat{M}}_k \\
\overline{\hat{\Phi}}_m
\end{bmatrix}
\end{equation}
where 
\begin{equation}
\overline{q}_h^d = \frac{1}{2}\sum_i
\Bigg(\alpha_i^{d;n}\hat{q}_i^{PG;n} + \alpha_i^{d;n+1}\hat{q}_i^{PG;n+1}\Bigg)
\end{equation}
is a time centered approximation to the upwinded potential vorticity, $q^d_h$, and
$\overline{\hat{M}}_h$, $\overline{\hat{\Phi}}_h$ are the degrees of freedom of the discrete 
mass flux and Bernoulli function exactly integrated between time levels $n$ and $n+1$ \cite{WCB20,Lee20}.
Note that the block matrix in \eqref{sw_mom_eq} is skew-symmetric, and so energy conservation is 
satisfied by the pre-multiplication of both sides by $[\overline{\hat{M}}_h^{\top}\quad\overline{\hat{\Phi}}_h^{\top}]$,
in analogy to the conservation of energy for the skew-symmetric advection equation \eqref{adv_eqn_ss}.

Since the potential vorticity, $q_h\in {W}_h$ is included within the skew-symmetric block matrix, 
the upwinding modification of the potential vorticity does not affect the conservation properties
of the model. While the upwinding of the potential vorticity does dissipate potential enstrophy,
this is not conserved in the original formulation due to the use of inexact quadrature 
\cite{LPG18}. This upwinded potential vorticity formulation is conceptually similar to the anticipated 
potential vorticity method \cite{SB85},
which similarly stabilises the potential vorticity by dissipating enstrophy. Note that this 
vorticity stabilisation does not address nonlinearities associated with the oscillation of 
gravity waves which may still yield aliasing errors, and also that the upwinding scheme
presented here is not strictly monotone, so some of these oscillations may still derive from the
rotational term in \eqref{sw_mom_eq}. In recent work \cite{WCB20} a method for determining 
a secondary velocity reconstruction within the discrete form of the skew-symmetric operator has been presented,
so as to upwind the potential gradients and mass fluxes associated with the generation of gravity
waves without dissipating energy.

\begin{figure}
\centering
\includegraphics[width=0.72\textwidth,height=0.54\textwidth]{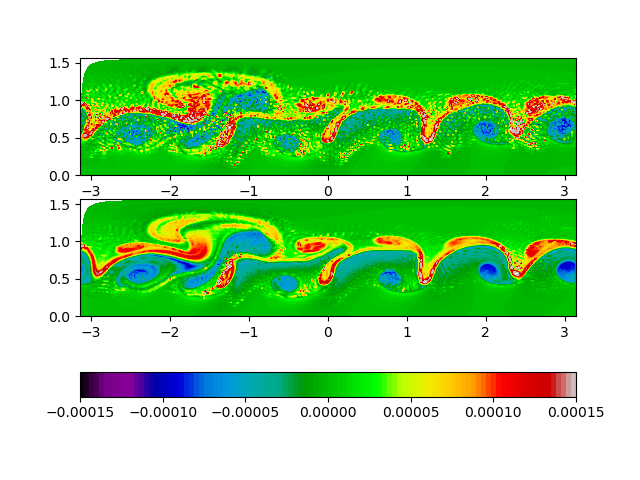}
\caption{Vorticity field for the Galewsky test case (day 7); top: original formulation, 
bottom: Petrov-Galerkin upwinded formulation. Only the northern hemisphere is shown.}
\label{fig::vorticity}
\end{figure}

The code was configured according to a standard test case for rotating shallow water on the sphere \cite{Galewsky04},
which is initialised as a steady jet overlaid with a small perturbation of the height field.
Over several days this acts to destabilise the jet into a series of barotropic eddies.
In each case the model was run with 32 elements of degree $p=3$ on each face of the cubed sphere
and a time step of $\Delta t = 120$s. 
As observed in Fig. \ref{fig::vorticity},
the Petrov-Galerkin upwinding of equation \eqref{vort_eqn_up} leads to a more coherent 
solution, with fewer aliasing errors.

\begin{figure}
\centering
\includegraphics[width=0.48\textwidth,height=0.36\textwidth]{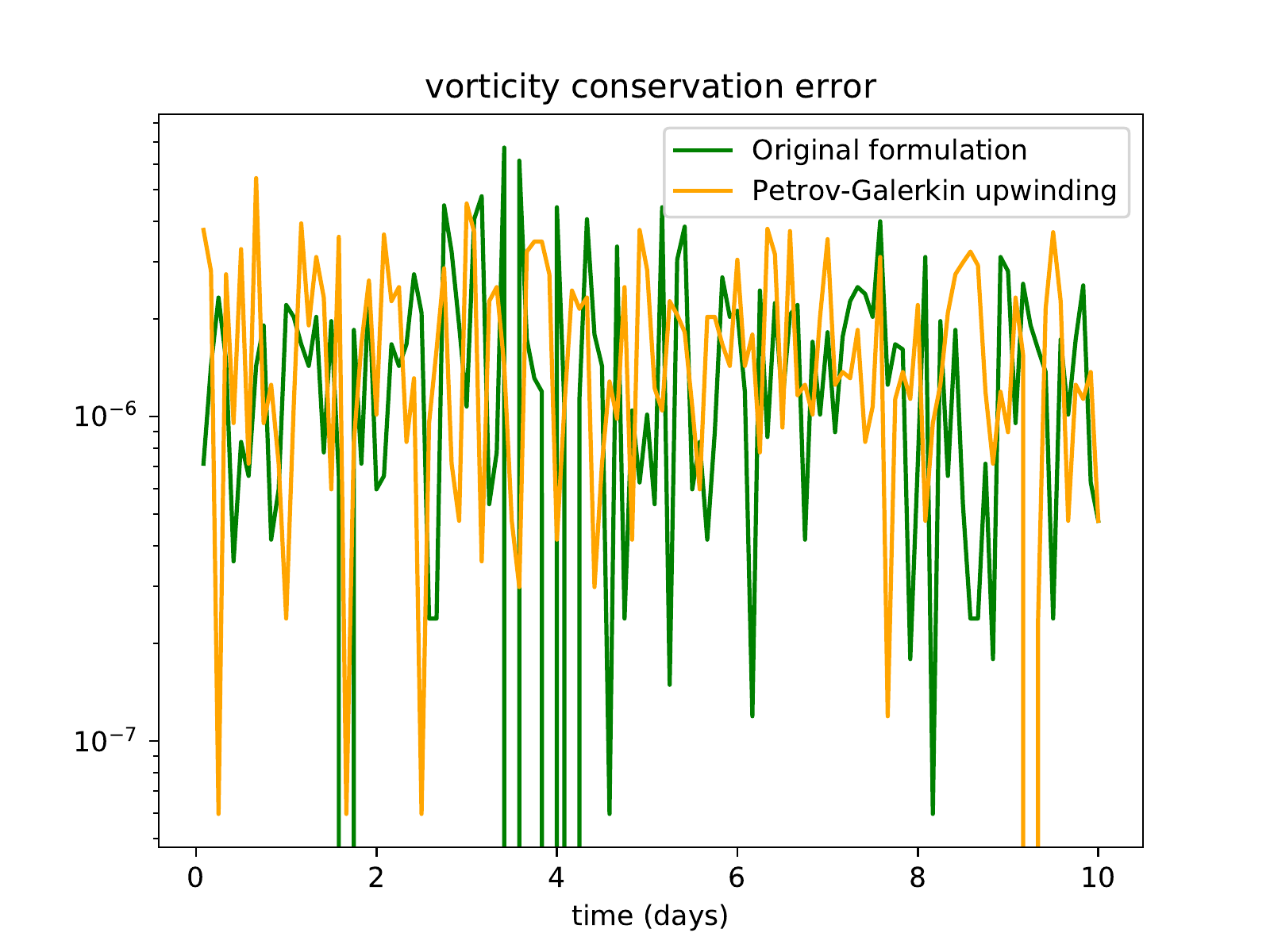}
\includegraphics[width=0.48\textwidth,height=0.36\textwidth]{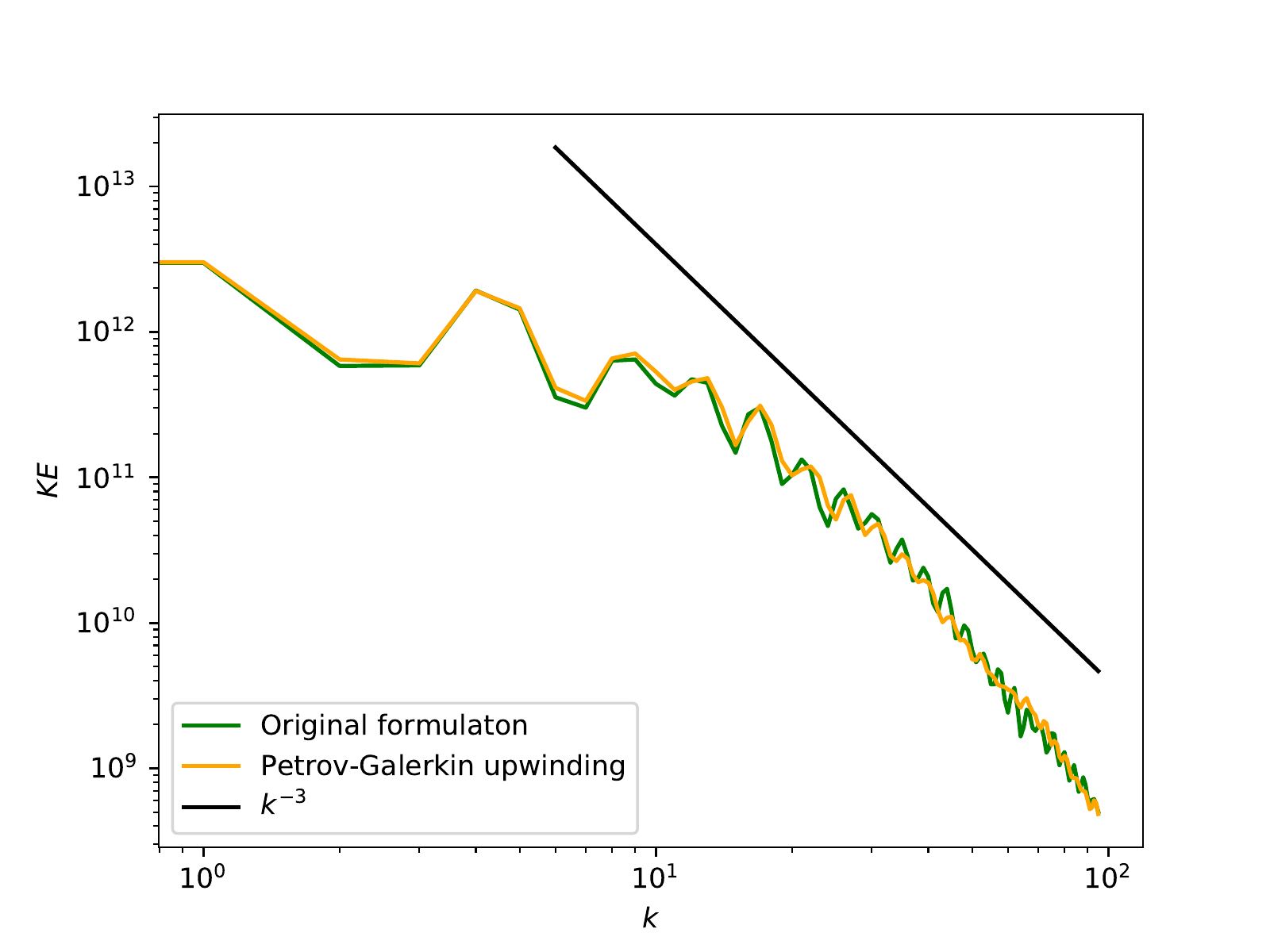}
\caption{Vorticity conservation errors (un-normalised) for the original and upwinded formulations (left), 
and kinetic energy power spectra at day 7 (right).}
\label{fig::vorticity_2}
\end{figure}

Figure \ref{fig::vorticity_2} shows the vorticity, $\omega$, conservation errors integrated over the sphere.
These are un-normalised, since the total vorticity integrates to zero, and integrated over
the sphere in physical units, such that truncation errors may be non-negligible.
There is no perceivable difference in the errors between the two schemes, and these are consistent
with previous results \cite{LP18}. This figure also shows the kinetic energy spectra of the two
schemes at day 7. These are computed by first interpolating the kinetic energy onto a regular
latitude--longitude grid, and then decomposing this solution into spherical harmonics. Both
schemes exhibit an upscale cascade of $k^{-3}$, consistent with the theory of two dimensional
turbulence, however the upwinded scheme shows less oscillation in this cascade, reflecting the minimisation 
of aliasing errors.

\section{Conclusions}

This article describes the formulation of upwinded advection operators for mixed mimetic spectral elements in
both flux form and material form. These operators exhibit dissipation of high wave numbers which suppresses 
spurious high frequency oscillations for sharp, poorly resolved gradients. The dissipation profiles increase 
with polynomial degree in a way that reflects the profile of higher order viscosity terms. 
Moreover the upwinded operators seal up the spectral gaps observed in dispersion relations for high 
order finite element methods. 
These upwinded formulations are relatively simple to compute and cheap to assemble. 

This upwinding formulation has been coupled to an existing solver for the shallow water equations on 
the sphere, where it is shown to reduce aliasing errors, as well as to improve the turbulence profile. In
future work the application of this scheme to the stabilisation of temperature fluxes for the compressible
Euler equations will be investigated.

The scripts used to generate the one dimensional results in Section 3 can be obtained from the author's
Github page at \verb|https://github.com/davelee2804/MiMSEM/tree/master/adv_eqn_1d|, while those for the two dimensional
test cases are found at \verb|https://github.com/davelee2804/MiMSEM/tree/hevi/sandbox/src|.

\section*{Appendix: computation of eigenvalues}

Assuming exact time integration, the semi-discrete advection equation is expressed as 
\begin{equation}
\boldsymbol{\mathsf M}^{-1}\boldsymbol{\mathsf A}q_h = \omega q_h
\end{equation}
for $q_h(x,t) = q_h(x)e^{-\omega t}$. The eigenvalues, $\omega_h$ and eigenvectors, $v_h$ of 
this square operator may be computed using any standard eigenvalue solver as 
$\omega_h,v_h = \tt{eig}(\boldsymbol{\mathsf M}^{-1}\boldsymbol{\mathsf A})$ (and
similarly for $\boldsymbol{\mathsf A}_{PG;\Delta t}$ and $-\boldsymbol{\mathsf A}_{PG;-\Delta t}^{\top}$).
We then define two additional operators, an interpolation operator between degrees of 
freedom in $\mathcal Q_h$ and physical coordinates, $\boldsymbol{\mathsf Q}$, for which
\begin{equation}
\mathsf Q_{jk} := e_k(\xi_j)
\end{equation}
(for which the Jacobian terms cancel), and a Fourier interpolation operator, $\boldsymbol{\mathsf F}$,
for which
\begin{equation}
\mathsf F_{jk} := \cos(2\pi k x_j/L) + i\sin(2\pi k x_j/L)
\end{equation}
where $L$ is the domain length, $i$ denotes the imaginary number, and $k$ a
given Fourier mode. Each eigenvector, $v_h$ is then projected onto a vector $v_h^f$ 
representing an expansion over Fourier modes as
\begin{equation}
v_h^f = \boldsymbol{\mathsf F}^{-1}\boldsymbol{\mathsf Q}v_h.
\end{equation}
We then sort $v_h^f$ in order to determine the Fourier mode with the largest
amplitude, $k_h$. This mode is then paired with the original eigenvalue, $\omega_h$ in order
to construct the dispersion relation.

\section{Acknowledgments}

David Lee would like to thank Dr. Darren Engwirda for the idea of applying the Petrov-Galerkin
scheme to the problem of potential vorticity advection in geophysical flows.
This project was supported by resources and expertise provided by CSIRO IMT Scientific Computing.

\end{document}